\input amstex.tex
\documentstyle{amsppt}
\UseAMSsymbols\TagsAsMath\widestnumber\key{ASSSS}
\magnification=\magstephalf\pagewidth{6.2in}\vsize8.0in
\parindent=6mm\parskip=3pt\baselineskip=16pt\tolerance=10000\hbadness=500
\NoRunningHeads\loadbold\NoBlackBoxes\nologo
\def\today{\ifcase\month\or January\or February\or March\or April\or May\or June\or
     July\or August\or September\or October\or November\or December\fi
     \space\number\day, \number\year}

\def\ns#1#2{\mskip -#1.#2mu}\def\operator#1{\ns07\operatorname{#1}\ns12}
\def\pa{\partial}\def\na{\nabla\!}
\def\dt{{\ns18{\Cal D}_{\ns23 t\!} }}\def\pdt{{\ns18\partial\ns18{\Cal D}_{\ns23 t\!} }}
\def\tr{\operatorname{tr}}
\def\div{\operator{div}}\def\curl{\operator{curl}}
\def\Vol{\operatorname{Vol}}\def\dist{\operatorname{dist} }
\def\leftalignspace{\mskip -45mu}\def\rightalignspace{\mskip -90mu}

\topmatter
\title
Well-Posedness for the Linearized Motion of an 
Incompressible Liquid with Free Surface Boundary
\endtitle
\author  Hans Lindblad 
\endauthor
\thanks  The author was supported in part by the National 
Science Foundation.  \endthanks
\address 
University of California at San Diego
\endaddress 
\date June 13, 2002\enddate 
\email
lindblad\@math.ucsd.edu
\endemail
\endtopmatter

\document 
\head{1. Introduction}\endhead 
We consider Euler's equations 
describing the motion of a perfect incompressible fluid in vacuum: 
$$\align 
\big(\partial_t +V^k\partial_{k} \big) v_j &=-\partial_{j} p,\quad j=1,...,n
\quad \text{in}\quad {\Cal D},\tag 1.1\\
\div V &=\partial_{k} V^k=0 \quad \text{in}\quad {\Cal D}\tag 1.2
\endalign 
$$
where $\partial_i=\pa/\pa x^i$ 
and ${\Cal D}=\cup_{\, 0\leq t\leq T}\, \{t\}\times{\Cal D}_t$,
${\Cal D}_t\subset{\Bbb R}^n$.
Here $V^k=\delta^{ki} v_i=v_k$
and we use the summation convention over repeated upper and lower indices.
The velocity vector field of the fluid is $V$, 
$p$ is the pressure and 
${\Cal D}_t$ is the domain the fluid occupies at time $t$. 
We also require boundary conditions on the free boundary 
$\partial {\Cal D}=\cup_{\, 0\leq t\leq T}\, \{t\}\times\partial {\Cal D}_t$; 
$$
\align p=0,\quad &\text{on}\quad \partial {\Cal D},\tag 1.3\\
       (\partial_t+V^k\partial_{k})|_{\partial {\Cal D}}
&\in T(\partial {\Cal D}),\tag 1.4
\endalign 
$$
Condition (1.3) says that the pressure $p$ vanishes outside the domain
and condition (1.4) says that the boundary moves with the velocity $V$ 
of the fluid particles at the boundary.

Given a domain 
${\Cal D}_0\subset \Bbb R^n$, that is homeomorphic to the unit ball,
 and initial data $v_0$, satisfying the constraint 
(1.2), we want to find a set 
${\Cal D}=\cup_{\, 0\leq t\leq T}\, \{t\}\times{\Cal D}_t$,
${\Cal D}_t\subset{\Bbb R}^n$
and a vector field $v$ solving (1.1)-(1.4) with initial conditions 
$$
 \{x;\, (0,x)\in {\Cal D} \}={\Cal D}_0,\qquad\text{and}\qquad 
v=v_0,\quad\text{on}\quad \{0\}\times {\Cal D}_0 \tag 1.5 
$$
Let ${\Cal N}$ be the exterior unit normal to the 
free surface $\pa{\Cal D}_t$. 
Christodoulou\cite{C2} conjectured
 that the initial value problem (1.1)-(1.5), 
is well posed in Sobolev spaces if
$$
\na_{\Cal N}\, p\leq -c_0 <0,\quad\text{on}\quad \partial {\Cal D},
\qquad \text{where}\quad \na_{\Cal N}={\Cal N}^i\partial_{x^i}. \tag 1.6
$$

Condition (1.6) is a natural {\it physical condition} since the pressure 
$p$ has to be positive in the interior of the fluid. 
It is essential for the well posedness in Sobolev spaces.
A condition related to Rayleigh-Taylor instability 
in \cite{BHL,W1} turns out to be equivalent to (1.6), see \cite{W2}. 
Taking the divergence of (1.1) gives:
$$
-\triangle p= (\partial_{j} V^k)\partial_{k} V^j,
\qquad \text{in}\quad {\Cal D}_t,\qquad p=0,\quad \text{on}\quad\partial
{\Cal D}_t\tag 1.7 
$$
In the irrotational case, when
$\curl {\,\ns18}v_{\ns10 \, ij}\!=\!\partial_i v_j\!-\partial_j v_i\!=\!0$, then 
 $\triangle p\!\leq\!0$ so $p\!\geq\!0$ and (1.6) holds by the strong
maximum principle. 
Wu \cite{W1,W2} proved well posedness locally in time,
(globally in space), in Sobolev spaces in the irrotational case.
Ebin \cite{E1} showed that the equations are ill posed when (1.6) is 
not satisfied and the pressure is negative and 
Ebin \cite{E2} announced an existence result when one adds 
surface tension to the boundary condition.
With Christodoulou $\ns40$\cite{CL}$\ns40$ we proved {\it $\ns30$a$\ns30$ priori$\ns10$} bounds in Sobolev
spaces in the general case of non vanishing curl, assuming (1.6). 
Usually if one has {\it a priori} estimates, existence follows from
similar estimates for some regularization or iteration scheme for the equation.  
However, the sharp estimates in \cite{CL} use all the symmetries of the 
equations and so only hold for perturbations of the equations that preserve the
symmetries. Here we show existence in Sobolev spaces for the linearized equations
using a new type of estimates.

The incompressible perfect fluid is to be thought of as an idealization of a liquid.
For small bodies like water drops surface tension should help holding it together and for 
larger denser bodies like stars its own gravity should play a role. 
Here we neglect the influence of such forces. Instead it is the incompressibility condition that prevents the 
body from expanding and it is the fact that the pressure is positive 
that prevents the body from breaking up in the interior.
Let us also point out that, from 
a physical point of view one can alternatively think of the pressure as being a small 
positive constant on the boundary instead of vanishing. 
What makes this problem difficult is that the regularity 
of the boundary enters to highest order. 
Roughly speaking, the velocity tells the boundary where to move and the boundary 
is the zero set of the pressure that determines the acceleration. 

Some existence results in Sobolev spaces are known in the irrotational case,
for the closely related water wave problem which describes
the motion of the surface of the ocean under the influence of earth's gravity.
In that problem, the gravitational field can be considered as uniform,  
however this problem reduces to our problem by going to an 
accelerated frame. The domain ${\Cal D}_t$ is unbounded for the water wave problem 
coinciding with a half-space in the case of still water.  
Nalimov\cite{Na} and Yosihara\cite{Y} proved local existence in Sobolev spaces 
in two space dimensions for initial conditions sufficiently close to still water.
Beale, Hou and Lowengrab\cite{BHL} have given an argument to show  
that problem is linearly well posed in a weak sense in Sobolev spaces,  
assuming a condition, which can be shown to be equivalent to (1.6). 
The condition (1.6) prevents the Rayleigh-Taylor instability 
from occurring when the water wave turns over. 
Finally Wu\cite{W1,2} proved local existence in general in two and three dimensions
for the water wave problem. The method of proofs in these papers uses 
that the velocity is irrotational and divergence free and hence harmonic 
to reduce the equations to equations on the boundary only.  

The main result here is existence for the linearized equations
in the case of non vanishing curl. 
The irrotational case was proved by Yosihara \cite{Y}.
The proof in \cite{Y}, see also \cite{W1,W2},
reduces the equation to the boundary and it does not generalize. 
Instead, we project the linearized equation onto an equation in the interior
using the orthogonal projection onto 
divergence free vector fields in the $L^2\!$ inner product.
This removes a difficult term, the differential 
of the linearization of the pressure, and reduces a higher order term, the
linearization of the moving boundary, to a symmetric unbounded
operator on divergence free vector fields.
The linearized equation becomes an evolution equation in the interior for this 
operator, which we call the normal operator.
It is basically the differential of the harmonic extension to the interior
of the normal component. In the irrotational case 
it becomes the normal derivative which is elliptic on harmonic functions 
and our equation reduces to an equation on the boundary similar to
those in \cite{Y,W1,W2}.

The normal operator is positive due to (1.6) and this will lead to energy bounds.
However, existence of regular solutions does not follow from standard energy 
methods or
semi-group methods since the operator is time dependent and non-elliptic in the case 
of non vanishing curl. Usually one gets equations and estimates for higher 
derivatives by commuting differential operators through the equation, 
but we can only use operators whose commutator with the 
normal operator is controlled by the normal operator. 
Geometric arguments lead us to use
Lie derivatives with respect to divergence free vector fields tangential 
at the boundary.
The commutators of these with the normal operator are controlled by the normal 
operator and they preserve the divergence free condition.
The same considerations apply to time differentiation so one should use the 
Lie derivative with respect to the material derivative (1.4) which reduces to
the time derivative of the vector field in the Lagrangian coordinates. 
To get estimates for all derivatives we 
use the fact that we have a better evolution equation for the curl and that any 
derivative can be controlled by tangential derivatives, the curl and the divergence.

As pointed out above, existence does not follow directly from estimates but one
must have existence and uniform estimates for some regularizing sequence.   
We replace the normal operator by a sequence of 
bounded operators converging to it which  are still 
symmetric, positive and they uniformly satisfy the same commutator estimates with 
the differential operators above.
Due to the geometric construction of the differential operators
there is a natural regularization which corresponds to
replacing the boundary by an inhomogeneous term supported in a 
small neighborhood of it.

Existence for the linearized equations or some modification
will be part of any existence
proof for the nonlinear problem.
The estimates here require more regularity of
the solution we linearize around than we get for the linearization.
However, we use the techniques presented here in a forthcoming paper \cite{L3},
to prove existence for the nonlinear problem with the Nash-Moser technique.

In order to formulate the linearized equations one has to introduce some
parametrization of the boundary. 
Let us therefore first express Euler's equations in the 
 Lagrangian coordinates in which the boundary becomes fixed. 
Given a domain ${\Cal D}_0$ in $\bold{R}^n$,
that is diffeomorphic to the unit ball $\Omega$, 
we can by a theorem in \cite{DM} find a
diffeomorphism $f_0:\Omega\to{\Cal D}_0$ that up to a constant factor
is volume preserving,  
i.e. after an additional scaling 
$\det(\partial f_0/\partial y)\!=\!1$. 
Assume that ${\Cal D}$ and $v\in C({\Cal D})$ are given satisfying (1.4).
The Lagrangian coordinates $x\!=\!x(t,y)\!=\!f_t(y)$ are given by solving 
$$
\frac{d x}{dt}=V(t,x(t,y)), \qquad x(0,y)=f_0(y),\qquad y \in \Omega.\tag 1.8
$$
Then $f_t:\Omega\to {\Cal D}_t$ is a volume preserving diffeomorphism,
since $\div V=0$,  and the boundary becomes fixed in the new $y$ coordinates. 
Let us introduce the notation 
$$
 D_t=\frac{\partial }{\partial t}\Big|_{x=const}\!\!
+\,  V^k\frac{\partial}{\partial x^k}
=\frac{\partial }{\partial t}\Big|_{y=const}\!\!
\qquad\text{and}\qquad \partial_i=\frac{\partial}{\partial {x^i}}
=\frac{\partial y^a}{\partial x^i}
\frac{\partial}{\partial y^a}, \tag 1.9
$$
for the material derivative and partial differential operators expressed in the 
Lagrangian coordinates. 
In  these coordinates Euler's equations (1.1),
the incompressibility condition (1.2) and the boundary condition (1.3) become
$$
D_t^2 x^i=-\pa_{i\,} p,\qquad \det{(\pa x/\pa y)}=1,
\qquad\text{in }\quad [0,T]\times\Omega\qquad\text{and}\qquad 
p\,\Big|_{\pa\Omega}=0.\tag 1.10
$$
where $p\!=\!p(t,y)$, $\pa_i$ now is to be thought of as the differential operator in 
(1.9) in $y$ and $D_t$ is the time derivative. We then define $V=D_t x$.
Note that the second equation in (1.10) follows since 
$D_t\ln{\big(\! \det{(\pa x/\pa y)}\big)}\!=
\div V\!=\!0$. 
Taking the divergence of the first equation in (1.${}_{\!}$10) gives (1.${}_{\!}$7) 
so $p$ is determined as functional of $(x,V)$. 
The initial conditions (1.5) become
$$
x\big|_{t=0}=f_0,\qquad D_t \, x\big|_{t=0}=V_0\tag 1.11
$$
subject to the constraints, 
$$
\det{(\pa f_0/\pa y)}=1,\qquad\text{and}\qquad
\div V_0=0\tag 1.12
$$
and Christodoulou's physical condition become
$$
\na_{\Cal N} p\big|_{\pa \Omega}\leq -c_0<0\tag 1.13
$$ 
where ${\Cal N}$ is the exterior unit normal to 
$\pa{\Cal D}_t$ parametrized by $x(t,y)$. 

Let us now derive the linearized equations of (1.10). 
We assume that $(x(t,y),p(t,y))$ is a given smooth solution of (1.10) satisfying (1.13)
for $0\leq t\leq T$. 
Let $\overline{x}(t,y,r)$ and $\overline{p}(t,y,r)$ be smooth functions 
also of a parameter $r$, such that  $(\overline{x},\overline{p})\big|_{r=0}=(x,p)$ and set
$(\delta x,\delta p)=(\pa \overline{x}/\pa r,\pa\overline{p}/\pa r)\big|_{r=0}$. 
Then the linearized equations is the requirement on $(\delta x,\delta p)$,
that $(\overline{x},\overline{p})$
satisfies the equations (1.10) up to terms bounded by $r^2$ as $r\to 0$. 
In other words, if 
$$
\Phi_i(x,p)=D_t^2 x^i+\pa_i p,\quad i=1,...,n,\qquad 
\Phi_0(x,p)=\det{(\pa x/\pa y)}-1,\qquad
\Phi_{n+1}(x,p)=p\big|_{\pa\Omega},\tag 1.14
$$
then then linearized operator is defined by
$$
\Phi^{\,\prime}(x,p)(\delta x,\delta p)
=\frac{\pa \Phi(\overline{x}, \overline{p})}{\pa r}\Big|_{r=0},
\qquad\text{where}\quad \overline{x}=x+r\delta x,\quad \overline{p}=p+r\delta p\tag 1.15
$$
Euler's equations (1.10) become $\Phi(x,p)=0$ and the linearized equations are
$$
 \Phi^{\,\prime}(x,p)(\delta x,\delta p)=0\tag 1.16
$$

Applying the operator $\delta f=\pa f/\pa r\big|_{r=0}$ to (1.10), 
using that by (2.8) $[\delta ,\pa_i]\!=\!-(\pa_i\delta x^k)\pa_k$,
gives the linearized equations 
$$
D_t^2 \delta x^i -(\pa_k p)\pa_i \delta x^k=-\pa_i \delta p ,
\qquad  \div\delta x=0,\qquad \text{and}\qquad\delta p\Big|_{\pa\Omega}=0.
\tag 1.17
$$
where we used that $\delta\ln{\big(\! \det{(\pa x/\pa y)}\big)}\!=
\div \delta x\!$, see (2.6). 
Here $\delta p$ is determined as a functional of $(\delta x,D_t\delta x)$ since taking  
the divergence of (1.17) gives an elliptic equation for $\delta p$ similar to (1.7). 
We now want to sow existence for (1.17) with
initial data 
$$
\delta x\big|_{t=0}=\delta f_0,\qquad
D_t\,\delta x\big|_{t=0}=\delta V_0,\tag 1.18
$$
satisfying the constraints
$$
\div \delta f_0=0,\qquad \div \delta V_0=(\pa_i \delta f_0^k)\pa_k V_0^i
\tag 1.19
$$
We remark, that the difference between (1.10) and (1.17) is the
term $\pa_k p\,\pa_i\delta x^k$ in (1.17). This term is higher order but
because of the sign condition (1.6) it will contribute with a positive term to the energy. 
We also remark that the equation (1.17) also shows up in estimating energies of higher 
order derivatives for (1.10) in \cite{CL}. In fact, the material derivative $D_t$ corresponds 
to the variation $\delta$ given by time translation. 
Our main result is: 
\proclaim{Theorem {1.}1} Let $\Omega$ be the unit ball in $\bold{R}^n$
and suppose that $(x,p)$ is a smooth solution of (1.10) satisfying
(1.13) for $0\leq t\leq T$. Suppose that $(\delta f_0,\delta V_0)$ are smooth
satisfying the constraints (1.19).
Then the linearized equations (1.17) have a smooth solution $(\delta x,\delta p)$
for $0\!\leq\! t\!\leq \!T$ satisfying the initial conditions (1.18). 
Let ${\Cal N}$ be the exterior unit normal to 
$\pa{\Cal D}_t$ parametrized by $x(t,y)$ and let  
$\delta x_{\Cal N}={\Cal N}\cdot\delta x$ be the normal component. Set
$$
E_r(t)=\| D_t \delta x(t,\cdot)\|_{H^r(\Omega)} +\|\delta x(t,\cdot)\|_{H^r(\Omega)}
+\|\delta x_{\Cal N}(t,\cdot)\|_{H^r(\pa \Omega)} \tag 1.20
$$
where
$H^r(\Omega)$ and $H^r(\pa\Omega)$ are the Sobolev spaces in $\Omega$ respectively on 
$\pa\Omega$. Then there are constants $C_r$ depending only on $(x,p)$, $r$ and $T$
such that 
$$
E_r(t)\leq C_r E_r(0),\qquad\text{for}\qquad 0\leq t\leq T ,\qquad r\geq 0. \tag 1.21
$$

Furthermore, let
$N^r(\Omega)$ be the completion of $C^\infty(\overline{\Omega})$ 
divergence free vector fields in the norm $\|\delta x\|_{H^r(\Omega)}
+\|\delta x_{\Cal N}\|_{H^r(\pa \Omega)} $. Then if the constraints in (1.19) hold and 
$$
(\delta f_0,\delta V_0)\in 
N^r(\Omega)\times H^r(\Omega)\tag 1.22
$$
 it follows that (1.17)-(1.18) has a solution 
$$
(\delta x, D_t\delta x)\in C([0,T],N^r(\Omega)\times H^r(\Omega)).\tag 1.23
$$
\endproclaim

As we have argued, any smooth solution of (1.1)-(1.5) 
with ${\Cal D}_0$ diffeomorphic to the unit ball can be reduced to a 
smooth solution of (1.10) where $\Omega$ is the unit ball. 
The term $\|\delta x_{\Cal N}\|_{H^r(\Omega)}$
is equivalent to the variation of the second fundamental form 
$\theta=\overline{\partial} {\Cal N}$ of the free boundary $\partial{\Cal D}_t$
measured in $H^{r-2}(\Omega)$, so our energy is essentially
$\|\delta\theta\|_{H^{r-2}(\pa\Omega)}+\|\delta v\|_{H^r(\Omega)}$.
This is to be compared with the {\it{ a priori}} bounds for the nonlinear
problem in \cite{CL} for 
$\|\theta\|_{H^{r-2}(\pa\Omega)}+\| v\|_{H^r(\Omega)}$. 
A slightly more general theorem holds, see section 2. 
Let us now outline the main ideas in the proof. 
We will rewrite the linearized equations (1.17) in a geometrically invariant way
and use this to obtain energy bounds and a regularization of the equation
which will give existence. 

We have defined our functions and vector fields to be functions of the Lagrangian
coordinates  $(t,y)\in[0,T]\times\Omega$ but we can alternatively think of them
as functions of the Eulerian coordinates
$(t,x)\in {\Cal D}$, and we will make this identification 
without explicitly saying that we compose with the inverse of the change of coordinate
$y\to x(t,y)$. The time derivative has a simple expression in the 
Lagrangian coordinates but the space derivatives have a simpler expression
in the Eulerian coordinates, see (1.9). For the most part we will think of our
functions and vector fields in the Lagrangian frame but we use the 
inner product coming from the Eulerian frame, i.e. in the Lagrangian frame
we use the pull-back metric of the Euclidean inner product:
$$
X\cdot Z=\delta_{ij} X^i Z^j=g_{ab} X^a Z^b,\qquad\text{where}\quad 
X^a=X^i\frac{\pa y^a}{\pa x^i},\qquad
 g_{ab}=\delta_{ij}\frac{\pa x^i}{\pa y^a}\frac{\pa x^j}{\pa y^b}.\tag 1.24
$$
Here $X^i$ refers to the components of the vector $X$ in the Eulerian frame,
$X^a$ refers to the components in the Lagrangian frame, $g_{ab}$ is the metric
in the Lagrangian frame and $\delta_{ij}$ is the Euclidean metric in the Eulerian frame.
The letters $a,b,c,d,e,f,g$ will refer to indices in the Lagrangian frame whereas
the indices $i,j,k,l,m,n$ will refer to the Eulerian frame. 
The norms and most of the operators we consider have an invariant interpretation
so it does not matter in which frame they are expressed. 
In the introduction we use express the vector fields in the Eulerian frame but
later we express the vector fields in the Lagrangian frame. 
The $L^2$ inner product of vector fields is given by 
$$\align 
\vspace{-0.0in}
\langle X,Z\rangle&=\int_{\dt} X\cdot Z \, dx=
\int_{\Omega} X\cdot Z\, dy\tag 1.25\\
\vspace{-0.09in}\endalign 
$$
where the equality follows from the incompressibility condition:
 $\det{(\pa x/\pa y)}=1$. 

We now want to derive energy bounds for the linearized equations (1.17).
Let us first point out that the boundary condition $p\big|_{\pa\Omega}=0$
implies that the energy is conserved for a solution of Euler's equations (1.10).
We have 
$$
\quad\frac{d}{dt} \int_\dt\!\! |V|^2  dx=\!\int_\dt\!\!\!  D_t |V|^2 dx
=-2\!\int_\dt\!\!\! V^i\pa_{i\,} p\, dx=2\!\int_\dt \!\!\div V\, p\, dx
-2\!\int_{\pdt} \!\!\!\! V_{ {\Cal N}\,} p\, dS=0,\!\!\!\!\!\tag 1.26
$$
where $V_{\Cal N}\!={\Cal N}_i V^i\!$ is the normal component of $V\!$. 
In fact, the first equality follows from the incompressibility condition
after expressing the integrals as
integrals over $\Omega$ as in (1.25),
 the second is Euler's equations (1.10), the third follows
from the divergence theorem and the last is the boundary condition and the 
divergence free condition. 

We will now use the orthogonal projection onto divergence free vector fields 
to rewrite the linearized equations (1.17) in an invariant way that can be 
used to derive energy bounds and for which there is a natural regularization. 
The orthogonal projection onto divergence free vector fields in the inner product 
(1.25) is given by 
$$
 PX^i=X^i -\delta^{ij}\pa_{\,j} {}_{\,} q ,\qquad \text{where}
\qquad \triangle q=\div X,\quad
q\Big|_{\pdt}=0. \tag 1.27
$$
We now want to project the first equation in (1.17) onto divergence free
vector fields. This removes the right hand side $\pa_i\delta p$, 
since we project along gradients of functions that vanish on the boundary. 
The second term in the first equation in (1.17) can be written as 
$-\pa_i \big( (\pa_k p)\delta x^k\big)+(\pa_i\pa_k p)\delta x^k$, where the 
last part is lower order and the projection of the first part turns out to
be a positive symmetric operator on divergence free vector fields. 
We define the normal operator $A$ to be 
$$
A X^i =P\big(-\delta^{ij}\pa_j (X^k \pa_k p)\big)
=-\delta^{ij}\pa_j\big(X^k \pa_k p-q\big).\tag 1.28
$$
where $q$ is chosen so that the divergence of $AX$ vanishes and $q$ vanishes on 
the boundary. 
Then $A$ is a positive symmetric operator on divergence free vector fields, 
if condition (1.6) holds. In fact,
if $X$ and $Z$ are divergence free then 
$$
\!\!\!\!\langle X,AZ\rangle=-\int_{\dt} \!X^i\pa_i (Z^k \pa_k p)\, dx=
\int_{\pdt}\!\! X_{\Cal N} Z_{\Cal N} \,(- \na_{\Cal N} p) \, dS, 
\quad\qquad X_{\Cal N}={\Cal N}_i X^i\! \tag 1.29
$$

There is one more issue we have to deal with before writing up the 
linearized equations (1.17) in a more pleasant form.
The time derivative $D_t$ does not preserve the divergence free
condition so we have to modify it so it does. The operator
$$
{\Cal L}_{D_t} X^i=D_t X^i-(\pa_k V^i)X^k
=\frac{\pa x^i}{\pa y^a}D_t  \Big(\frac{\pa y^a}{\pa x^k} X^k\Big)\tag 1.30
$$
preserves the divergence free condition if $V$ is divergence free.
This is because it is the space time Lie derivative with respect to the divergence free
vector field $D_t=(1,V)$ restricted to the space components. Another way to look at it is
that it is just the time derivative of the vector field $X$ expressed 
in the Lagrangian frame. The divergence 
is invariant under coordinate changes and the volume form
is time independent so it commutes with time differentiation in the Lagrangian
coordinates. 

We now project the linearized equations (1.17) and get 
an evolution equation on divergence free vector fields for the normal operator $A$: 
$$
\ddot{X}^i+ A X^i=
 -2 P\big( (\pa_k V^i)\dot{X}^k \big) 
\qquad\text{where}\qquad  X=\delta x,\quad \dot{X}=
{\Cal L}_{D_t} \delta x, \quad 
\ddot{X}={\Cal L}_{D_t}^2 \delta x
\tag 1.31
$$
Introducing the orthogonal projection onto divergence free
vector fields solved to problems. First it turned the higher order term,
the second term in (1.17) into a positive symmetric operator. 
Secondly it got rid of the third term in (1.17) which caused considerable 
difficulties in \cite{CL}. In fact, the projection of a gradient of a function 
that vanishes on the boundary vanishes. The right hand side of (1.31) is lower 
order since the projection is a bounded operator. 
Associated with (1.31) is the energy
$$
E(t)=\langle \dot{X}, \dot{X} \rangle
+\langle X, (A+I) X \rangle\tag 1.32
$$
and one can show an energy estimate $|E^{\,\prime}(t)|\leq CE(t)$ which 
gives an energy bound. We remark that for divergence free vector fields 
(1.32) is equivalent to (1.20) with $r=0$. 
In order to show this energy bound we must calculate 
the commutator of the time derivative and the normal operator, which follows
from the argument below. 

In order to prove the energy bound and similar energy bounds for higher
derivatives one has to control the commutator of differential operators with the
normal operator. This is however a delicate matter since these commutators have 
to be controlled by the normal operator itself and only certain 
geometric operators satisfy this. Let $T$ be a divergence free vector field 
that is tangential at the boundary and let 
$$
{\Cal L}_T X^i=T^k\pa_k X^i-X^k\pa_k T^i\tag 1.33
$$ 
be the Lie derivative with respect to $T$ applied to a vector field $X$.
Then ${\Cal L}_T X$ is divergence free if $X$ is divergence free.  
It turns out that the commutators between ${\Cal L}_T$ and the normal operator
can be controlled by the normal operator:
$$
[{\Cal L}_T, A] X^i=( {\Cal L}_T \delta^{ij})\delta_{jk} A X^k
+A_{Tp} X^i\tag 1.34
$$
where for $f$ vanishing on the boundary we defined 
$$
A_f X=-P\big(\delta^{ij}\pa_j ( X^k \pa_k f)\big).\tag 1.35
$$
(1.34) follows from (1.28) using that the Lie derivative commutes with 
exterior differentiation and that the tangential 
derivatives $Tp$ and $Tq$ also vanish on the boundary since $p$ and $q$ do. 
 In view of the physical condition (1.13)
it follows from (1.29) that 
$$
|\langle X, A_{f} X\rangle| \leq C\langle X,A X\rangle,\qquad\quad\text{where}\qquad
C=\|\na_{\Cal N} f/\na_{\Cal N} p\|_{L^\infty(\pa\Omega)}.\tag 1.36
$$

Applying ${\Cal L}_T$ to the linearized equations (1.31) therefore
gives a similar equation for ${\Cal L}_T \delta x$ for which we also 
get energy bounds if $T$ is a divergence free 
vector field that is tangential at the boundary. 
The second term in the commutator (1.34) can be controlled using (1.36).
In order to control the first term in the commutator one has to use that 
$AX$ can be controlled in terms of ${\Cal L}_{D_t}^2 \delta x$ through the equation 
(1.31). Therefore we also have to differentiate the equation with respect to time
and include time derivatives up to highest order in the energies. 
We define energies
$$
E_r^{\Cal T}(t)= \sum_{|I|\leq r,\, I\in {\Cal T}}\sqrt{
\langle {\Cal L}_{D_t}{\Cal L}_T^I X, {\Cal L}_{D_t}{\Cal L}_T^I X
\rangle + \langle {\Cal L}_T^I X, A {\Cal L}_T^I X
\rangle },\qquad X=\delta x\tag 1.37 
$$ 
where ${\Cal T}$ is a family of divergence free vector fields that are tangential
at the boundary and span the tangent space of the boundary including the time 
derivative $D_t$ and ${\Cal L}_T^I$ is any product of $r=|I|$ Lie derivatives
with respect to these. Then one can prove energy estimates
$
E_r^{\Cal T}(t)\leq C E_r^{\Cal T}(0).
$

\comment 
$$
\underline{A} X_i=\delta_{ij} A X^j=-\pa_i\big( (\pa_k p) X^k\big)-\pa_i q\tag 1.35
$$
where by (1.28) 
$q$ is chosen so $q$ vanishes on the boundary and the divergence of $AX$ vanishes.
Then the commutator satisfies 
$$
P {\Cal L}_T\underline{A} X-\underline{A}{\Cal L}_T X=
A_{Tp} X\tag 1.33
$$
\endcomment 

The energies (1.37) only contain tangential derivatives.
In order to control normal derivatives also we use: 
$$
|\pa Z|\leq C\big(|\div Z|+|\curl \underline{Z}|+\sum_{S\in {\Cal S}} |SZ|\big)\tag 1.38
$$
where ${\Cal S}$ is a family of vector fields that span the tangent space
of the boundary and 
$\curl \underline{Z}_{ij}=\pa_i \underline{Z}_j-\pa_j\underline{Z}_i$, where 
$\underline{Z}_i=\delta_{ij} Z^j$ is the one form corresponding to the vector field $Z$.
 The divergence of $\dot{X}={\Cal L}_{D_t}\delta x$ 
vanishes and there is a better
evolution equation for $\curl \underline{\dot{X}}$.
In fact the curl of the higher order operator
$A$ in (1.31) considered as an operator with values in the one forms
vanishes since it is a gradient. For a solution of Euler's equations (1.10) the 
curl is preserved: 
$$
{\Cal L}_{D_t} \curl v=0,\tag 1.39
$$
where ${\Cal L}_{D_t}$ is the space time Lie derivative with respect to 
$D_t\!=\!(1,\!V)$ of the two form $\sigma$: 
$$
{\Cal L}_{D_t} \sigma_{ij}=D_t \,\sigma_{ij}
+(\pa_i V^l) \sigma_{lj}+(\pa_j V^l)\sigma_{il}
=\frac{\pa y^a}{\pa x^i}\frac{\pa y^b}{\pa x^j}
D_t \Big(\frac{\pa x^k}{\pa y^a}\frac{\pa x^l}{\pa y^b} \sigma_{kl}\Big),\tag 1.40
$$
restricted to the space components, i.e. it is the time derivative of the two 
form expressed in the Lagrangian frame. For the linearized equations we have the
following identity:
$$
{\Cal L}_{D_t}\curl\,\delta z=0,\qquad \qquad 
\delta z_i =\delta_{ij}{\Cal L}_{D_t} X^j-\curl v_{ij}\, X^j,
\qquad X^i=\delta x^i
\tag 1.41 
$$
Since the Lie derivative commutes with exterior differentiation 
$\curl {\Cal L}_T^I \delta z$ is also conserved. 

The above argument gives energy bounds, assuming existence.
However, existence does not follow directly from estimates.
To show existence we must approximate the
linearized equations with some equation for which we know there is existence 
and prove that we have uniform bounds for the norms as the approximation 
gets better so that we can construct a sequence that
tends to a solution of the linearized equations.
For $f>0$ in ${\Cal D}_t$ and $f\big|_{\pa {\Cal D}_t}=0$ we 
 define the smoothed out normal operator by
$$
A_f^\varepsilon X^i =P\big(- \chi_\varepsilon (d)\delta^{ij} \pa_j
 ( fd^{-1}X^k\pa_k d)\big)=
P\big(\chi_\varepsilon^\prime (d)\delta^{ij} (\pa_j d)
  fd^{-1}X^k\pa_k d\big)\tag 1.42
$$
where $d=d(y)=\dist{(y,\pa\Omega)}$ and $\chi_\varepsilon(d)=\chi(d/\varepsilon)$.
Here $\chi$ is a smooth cut off function,  $\chi(s)=1$, when $s\geq 1$,
$\chi(s)=0$, when $s\leq 0$ and $\chi^\prime(s)\geq 0$. 
Then $A_f^\varepsilon $ is a positive symmetric operator on divergence free 
vector fields, if condition (1.6) holds. In fact,
if $X$ and $Z$ are divergence free then 
$$
\!\!\!\!\langle X, A_f^\varepsilon Z\rangle=-\int_{\dt} \!\chi_\varepsilon(d)
X^i\pa_i (f d^{-1} Z^k \pa_k d)\, dx=
\int_{\dt}\!\! (X^i\pa_i d) (Z^k\pa_k d) \, \chi_\varepsilon^{\,\prime}
(d) f d^{-1} \, dx.\tag 1.43
$$
It follows that $A_f^\varepsilon$ is symmetric and positive and satisfies the
same commutator properties as $A_f$ and the curl of $A_f^\varepsilon$
vanishes when $d\geq \varepsilon$. Furthermore $A_f^\varepsilon $ is a bounded operator,
i.e. $\| A_f^\varepsilon X\|_r\leq C_{\varepsilon r} \|X\|_r$.

We will actually first obtain energy estimates for the linearized equations 
with vanishing initial data and an inhomogeneous divergence free
term that vanishes to any
order as $t\to 0$:
$$
\ddot{X}^i+ A X^i
 +2 P\big( (\pa_k V^i)\dot{X}^k \big) =\delta \Phi,
\qquad\quad {\Cal L}_{D_t}^k X\big|_{t=0}=0,\quad k\leq r,
\qquad \div\delta \Phi=0,
\tag 1.44
$$
of the form 
$$
E_r^{\Cal T}(t)\leq C_r\int_0^t\|\delta\Phi\|_{r}^{\Cal T}\,
d\tau ,\qquad\quad\text{where}\quad \|\delta\Phi\|_{r}^{\Cal T}
=\sum_{|I|\leq r,\, I\in {\Cal T}}\|{\Cal L}_T^I \delta\Phi\|\tag 1.45
$$
One can reduce to this situation by subtracting a power series solution 
in time to (1.31). 
 (1.44) with 
$A$ replaced by $A^\varepsilon=A_p^\varepsilon$ is just an ordinary differential equation
in $H^r(\Omega)$ so existence for this equation follows. Because $A^\varepsilon$
uniformly satisfies the same commutator estimates as $A$ we will obtain uniform 
energy bounds and will be able to pass to the limit as $\varepsilon\to 0$ and 
obtain a solution for (1.44). The reason we have to first subtract off the initial 
conditions in this way is that
the energy (1.37) contains time derivatives up to highest order and 
these would have to be obtained from the equation. 
The operator $A^\varepsilon$ is smoothing but only in the tangential directions
and in the normal directions it is worse than $A$ so if we had replaced $A$ by 
$A^\varepsilon$ directly in (1.31) the higher order initial conditions would 
have depended on $A^\varepsilon$ in an uncontrollable way. 
As described above, we will first prove the energy bounds 
 in such a way that we can obtain the same uniform
bounds for the smoothed out equation and pass to the limit as $\varepsilon\to 0$
to obtain existence. Once we have existence we can then obtain the more natural 
energy bounds for the initial value problem in Theorem {1.}1. 

\comment
First we introduce the Lagrangian coordinates and the
linearized equations in section 2, then we introduce the projection and the normal operator
in section 3
and after that prove the lowest order basic energy estimate in section 4. 
As it turns out, in order to prove existence for the initial value problem
we have to turn the initial conditions into an inhomogeneous term, in section 5.
Then in section 6 we construct the tangential vector fields, in section 7 we
calculate commutators between Lie derivatives with respect to tangential vector fields.
\endcomment

\head 2. Lagrangian coordinates, the linearized equation and statement of the theorem.\endhead 
Let us introduce Lagrangian coordinates in which the boundary becomes fixed.
Let $\Omega$ be a domain in $\bold{R}^n$ and let
$f_0:\Omega\to {\Cal D}_0$ be a diffeomorphism that is volume preserving;
$\det(\partial f_0/\partial y)=1$. For simplicity we will assume that 
$\Vol({\Cal D}_0)$ is the volume of the unit ball in $R^n$.
 By a theorem of \cite{DM} we can prescribe the volume form 
up to a constant for any mapping of one domain into another so
we may assume that $\Omega$ is the unit ball. 
Assume that
$v(t,x)$ and $p(t,x)$, $(t,x)\in  {\Cal D}$ are given satisfying the boundary
conditions (1.3)-(1.4). 
The Lagrangian coordinates $x=x(t,y)=f_t(y)$ are given by solving 
$$
{d x}/{dt}=V(t,x(t,y)), \qquad x(0,y)=f_0(y),\quad y \in \Omega\tag 2.1
$$
Then $f_t:\Omega\to {\Cal D}_t$ is a volume preserving diffeomorphism,
since $\div V=0$,  and the boundary becomes fixed in the new $y$ coordinates. 
Let us introduce the notation 
  $$\align
\vspace{-0.06in} 
\quad D_t&=\frac{\partial }{\partial t}\Big|_{y=constant}=
\frac{\partial }{\partial t}\Big|_{x=constant}
+\,  V^k\frac{\partial}{\partial x^k},\tag 2.2\\
\vspace{-0.3in}
&
\endalign 
$$
for the material derivative and 
$$\align 
\vspace{-0.07in}
\partial_i&=\frac{\partial}{\partial {x^i}}=\frac{\partial y^a}{\partial x^i}
\frac{\partial}{\partial y^a}. \tag 2.3
\endalign 
$$
In  these coordinates Euler's equation (1.1),
the incompressibility condition (1.2) and the boundary condition (1.3) become
$$
D_t^2 x^i=-\pa_i p,\qquad \quad \kappa=\det{(\pa x/\pa y)}=1,
\qquad\qquad   p\big|_{\pa\Omega}=0\tag 2.4
$$
where $x\!=\!x(t,y)$, $p\!=\!p(t,y)$.
The initial conditions (1.5) become
$$
 x\,\big|_{t=0}=f_0,\qquad
 D_t \, x\,\big|_{t=0}=V_0.\tag 2.5
$$
In fact, recall that $D_t\det{(M)}=\det{(M)}\tr\,( M^{-1} D_t{M})$,
for any matrix $M$ depending on $t$ so 
$$
 D_t  \det{(\pa x/\pa y)}= \det{(\pa x/\pa y)}\,
({\pa y^a}/{\pa x^i})({\pa D_t x^i}/{\pa y^a})=\pa_i D_t x^i=\div D_t x=\div
V=0.\tag 2.6 
$$
Note that $p$ is uniquely determined as a functional of $x$ by (2.4)-(2.5). 
In fact taking the divergence of
Euler's equations (2.4) using (2.6) 
gives $\triangle p=-(\pa_i D_t x^j)(\pa_j D_t x^i)$. 

Let $\delta$ be a {\it variation} with respect to some parameter $r$,
in the Lagrangian coordinates:
$$
\delta ={\partial}/{\partial r}\big|_{(t,y)=const},\tag 2.7
$$
We think of $x(t,y,r)$ and $p(t,y,r)$ as depending on 
$r$ and differentiate with respect to $r$. 
Differentiating (2.3) using the formula for the derivative of 
the inverse of a matrix, $\delta M^{-1}=-M^{-1}(\delta M)M^{-1}$, gives
$$
[\delta,\pa_i]=-(\pa_i \delta x^k)\pa_k.\tag 2.8
$$
Differentiating (2.4), using (2.8) and (2.6) with $D_t$ replaced by $\delta$
gives the linearized equations: 
$$
D_t^2 \delta x^i-(\pa_k p) \pa_i \delta x^k=-\pa_i \delta p,\qquad\quad
 \div \delta x=0,
\qquad\quad \delta  p\big|_{\pa\Omega}=0.\tag 2.9
$$ 
It is however better to use the fact that $v$ and $p$ are solutions of
Euler's equations, $D_t v_i=-\pa_i p$, to arrive at the following equation 
$$
D_t^2 \delta x^i-\pa_i \big((\pa_k p) \delta x^k\big) 
=-\pa_i \delta p+(\pa_k D_t v_i ) \delta x^k,
\qquad\quad \div \delta x=0,\qquad\quad \delta  p\big|_{\pa\Omega}=0.\tag 2.10
$$

We will now transform the vector field $\delta x$
to Lagrangian coordinates, because in these coordinates the time derivative
preserves the divergence free condition. Let 
$$
W^a=\delta x^i\, \frac{\pa y^a}{\pa x^i},\qquad 
\delta x^i= W^b\frac{\pa x^i}{\pa y^b},
\qquad q=\delta p.\tag 2.11
$$
The letters $a,b,c,d,e,f$ will refer to quantities in the Lagrangian frame whereas
the letters $i,j,k,l,m,n$ will refer to ones in Eulerian frame, e.g. 
$\pa_a=\pa/\pa y^a$ and $\pa_i=\pa/\pa x^i$.
With this convention we have 
$$
\pa_i
=\frac{\pa y^a}{\pa x^i}\pa_a ,\qquad
\pa_a =\frac{\pa x^i}{\pa y^a}\pa_i .\tag 2.12
$$
Multiplying the first equation in (2.10) by $\pa x^i/\pa y^a$ and summing over $i$ gives
$$
\delta_{ij}\frac{\pa x^i}{\pa y^a} D_t^2 \delta x^j
-\pa_a \big( (\pa_c p) W^c\big) =\pa_a q+\frac{\pa x^i}{\pa y^a} (\pa_c D_t v_i)W^c\tag 2.13
$$
since $(\pa_k p) \delta x^k\!=\!(\pa_c p) W^c$ and 
$(\pa_k D_t v_i ) \delta x^k\!=\!(\pa_c D_t v_i)W^c$. 
On the other hand 
$$\align 
D_t \,\delta x^i&= (D_t W^b)\frac{\pa x^i}{\pa y^b}+ 
 W^b \frac{\pa V^i}{\pa y^b},\qquad\qquad\text{and}  \tag 2.14\\
D_t^2 \delta x^i&= (D_t^2 W^b)\frac{\pa x^i}{\pa y^b}+
2\frac{\pa V^i}{\pa y^b} D_t W^b + W^b \frac{\pa D_t V^i}{\pa y^b}\tag 2.15
\endalign 
$$
Multiplying (2.15) by $\pa x^i/\pa y^a$, summing over $i$,
 and substituting into (2.13) gives
$$
\delta_{ij}\frac{\pa x^i}{\pa y^a}\frac{\pa x^j}{\pa y^b} D_t^2 W^b
-\pa_a \big( (\pa_c p) W^c\big)=\pa_a q-2\frac{\pa x^i}{\pa y^a}
\frac{\pa x^k}{\pa y^b} (\pa_k v_i )D_t W^b, \tag 2.16
$$
where
$$
g_{ab}=\delta_{ij}\frac{\pa x^i}{\pa y^a}\frac{\pa x^j}{\pa y^b} \tag 2.17
$$
is the metric $\delta_{ij}$ expressed in the Lagrangian coordinates. Let $g^{ab}$ be the inverse 
of the metric $g_{ab}$, 
$$
\dot{{g}}_{ab}= D_t g_{ab}=\frac{\pa x^i}{\pa y^a}
\frac{\pa x^k}{\pa y^b} \big(\pa_k v_i+\pa_i v_k\big)\qquad\text{and}\qquad 
\omega_{ab}= \frac{\pa x^i}{\pa y^a}
\frac{\pa x^k}{\pa y^b} (\pa_i v_k-\pa_k v_i)\tag 2.18 
$$
be the time derivative of the metric and the vorticity
in the Lagrangian coordinates. Expression (2.16) becomes
$$
g_{ab} D_t^2 W^b-\pa_a\big( (\pa_c p) W^c\big) 
=-\pa_a q-\big(\dot{{g}}_{ac}-\omega_{ac} \big) D_t W^c.\tag 2.19
$$

(2.19) can alternatively be expressed, using the inverse $g^{ab}$ of $g_{ab}$,
in the form 
$$
D_t^2 W^a-g^{ab}\pa_b\big( (\pa_c p) W^c\big)
=-g^{ab}\pa_b q-g^{ab}\big(\dot{{g}}_{bc}-\omega_{bc} \big) D_t W^c.\tag 2.20
$$
The divergence is invariant under coordinate changes so 
the second condition in (2.10) is 
$$
\div W=\kappa^{-1}\pa_a (\kappa W^a )=0,\qquad\text{where}\qquad
 \kappa=\det{(\pa x/\pa y)}=1\tag 2.21
$$
Finally, the last equation in (2.10) is, since $q=\delta p$,
$$
q\big|_{\pa\Omega}=0\tag 2.22
$$
Then linearized equations are now the requirement that (2.20), (2.21) and (2.22) hold
and we want to find $(W,q)$ satisfying these equations and the initial conditions
$$
W\big|_{t=0}=W_0,\qquad \dot{W}\big|_{t=0}=W_1,
\qquad\text{where}\qquad \div W_0=\div W_1=0,\qquad \dot{W}=D_t W
\tag 2.23
$$

We can however express (2.20)-(2.23) in as one equation as follows. 
First we note that $q=\delta p$ is determined as a functional of $W$ and $D_t W$. 
In fact, it follows from (2.21) that 
$\div D_t^2 W=0$ so taking the divergence of (2.20) using (2.22)
gives us an elliptic equation for $q$: 
$$
\triangle q=\kappa^{-1}\pa_a\big(\kappa  g^{ab}\pa_b q\big)=
\kappa^{-1}\pa_a\Big(\kappa g^{ab}\pa_b\big( (\pa_c p) W^c\big)
 -\kappa g^{ab}\big(\dot{{g}}_{bc}-\omega_{bc} \big) D_t W^c\Big),
\qquad  q\big|_{\pa \Omega}=0. \tag 2.24
$$
We now write $q=q_1+q_2+q_3$, where $q_i\big|_{\pa\Omega}=0$
and $\triangle q_i$ is equal to each of the three terms in the right 
hand side of (2.24). The equations (2.20)-(2.22) can then be written as one
equation, $L_1 W=0$ where 
$$
L_1 W=\ddot{W}+A W+\dot{G}\dot{W}-C\dot{W}\tag 2.25
$$
and
$$\align 
AW^a&=-g^{ab}\pa_b\big( (\pa_c p) W^c-q_1\big),
\qquad\quad \triangle q_1=\triangle \big((\pa_c p) W^c\big).
\qquad\quad q_1\big|_{\pa\Omega}=0\rightalignspace\tag 2.26\\
\dot{G} \dot{W}^a &= g^{ab}(\dot{g}_{bc} \dot{W}^c+q_2),\qquad\qquad\qquad
\triangle q_2= -\pa_a\big(g^{ab}\dot{g}_{bc} \dot{W}^c\big)
\qquad q_2\big|_{\pa\Omega}=0\rightalignspace\tag 2.27\\
C \dot{W}^a &=g^{ab}(\omega_{bc} \dot{W}^c-q_3),\qquad\quad\qquad\quad
\triangle q_3= \pa_a\big(g^{ab}\omega_{bc} \dot{W}^c\big)
\qquad\quad q_3\big|_{\pa\Omega}=0\rightalignspace\tag 2.28
\endalign
$$

We will prove the following theorem:
\proclaim{Theorem {2.}1}  Suppose that
$x, p\in C^{\infty}([0,T]\times\overline{\Omega})$, $p\, \big|_{\pa\Omega}=0$,
$\na_N p\, \big|_{\pa\Omega}\leq -c_0<0$ and $\div D_t x=0$. 
Suppose that $F\in C^{\infty}([0,T]\times\overline{\Omega})$ and $W_0, W_1
\in C^{\infty}(\overline{\Omega})$ are all divergence free. 
Then
$$
L_1 W=F,\qquad\quad W\big|_{t=0}=W_0,\qquad \dot{W}\big|_{t=0}=W_1,\tag 2.29
$$
where $L_1$ be given by (2.25)-(2.28), 
has a divergence free solution $W\in C^{\infty}([0,T]\times\overline{\Omega})$. 

Let $H^r(\Omega)$ be the Sobolev spaces and let
$N^r(\overline{\Omega})$ be the completion of 
$C^\infty(\overline{\Omega})$ divergence free vector fields in the 
norm $\|W\|_{H^r(\Omega)} +\|W_N\|_{H^r(\pa\Omega)}$, where 
$W_N=W\cdot N$ is the normal component. Then if 
$$
(W_0,W_1)\in N^r(\Omega)\times  H^r(\Omega),
\qquad F\in L^1\big([0,T],H^r(\Omega)\big) \tag 2.30
$$
are all divergence free 
it follows that (2.29) have a a divergence free solution
$$
(W,\dot{W})\in C\big([0,T],
 N^r(\Omega)\times H^r(\Omega)\big).\tag 2.31
$$
Moreover, with a constant $C$ depending only on the $C^{r+2}$ norm of $x$ and $p$ 
and the constant $c_0$ we have 
$$
\|\dot{W}(t)\|_{H^r} +\|W(t)\|_{N^r}
\leq C\Big( \|\dot{W}(0)\|_{H^r} +\|W(0)\|_{N^r}+\int_0^t\|F(\tau)\|_{H^r}\, d\tau\Big).
\tag 2.32
$$
\endproclaim
\demo{Remark} The restrictions that $\div V=0$ and $\div F=0$ can be 
removed and in order to use the Nash-Moser technique one indeed needs to show that
the linearized operator is invertible away from a solution and outside the divergence 
free class. In \cite{L3} the techniques presented here are used to show this. 
\enddemo

\head 3. The projection onto divergence free vector fields and 
the normal operator. \endhead 
Let $P$ be the orthogonal projection onto divergence free vector fields 
in the inner product 
$$
\langle W,U\rangle =\int_\Omega g_{ab} W^a U^b\, dy, \tag 3.1
$$
Then the projection $P$ 
$$
P U^a= U^a-g^{ab}\pa_b\,  q,\qquad 
\triangle q= \pa_a \big( g^{ab} \pa_b \, q\big)=\div U
= \pa_a ( U^a),\qquad q\Big|_{\pa\Omega}=0. \tag 3.2
$$
That this is the orthogonal projection follows since $g_{ab}g^{bc}=\delta_{a}^{c}$ 
and 
$$
\langle W,(I-P) U\rangle =-\int_\Omega g_{ab} W^a  g^{bc}\pa_c q\, dy
=\int_\Omega (\pa_a  W^a )\, q\, dy
-\!\!\int_{\pa\Omega} \!\!\! N_a W^a \, q \, dS=0,\quad \text{if}\quad \pa_a  W^a =0
\tag 3.3 
$$
where $N_a$ is the exterior unit conormal and $dS$ is the surface measure.
The projection of a gradient of a function that vanishes on the boundary vanishes:
$$
P\big( g^{ab}\pa_b q\big)=0,\qquad\quad\text{if}\qquad q\big|_{\pa\Omega}=0.
\tag 3.4 
$$

The projection has norm one:
$$
\| P U\|\leq  \|U\|,\qquad\quad
\|(I-P)U\|\leq \|U\|,
\qquad\quad \|W\|=\langle W,W\rangle^{1/2}\tag 3.5
$$
The projection is continuous on the Sobolev spaces $H^r(\Omega)$ if the metric is
sufficiently regular:
$$
 \| P U\|_{H^r(\Omega)}\leq C_r\|   U\|_{H^r(\Omega)},
\tag 3.6
$$
since it is just a matter of solving the Dirichlet problem: 
$$
\|q\|_{H^{r+1}(\Omega)}\leq C_r \|U\|_{H^r(\Omega)},\qquad r\geq 0,
\qquad\text{if}\quad \triangle q=\div U,\quad q\big|_{\pa\Omega}=0.\tag 3.7
$$
For $r\geq 1$ this is the standard estimate for the Dirichlet problem. 
For $r=0$ this is obtained by multiplying by $q$, 
using that the right hand side is in divergence form, 
integrating by parts and using that $q\big|_{\pa\Omega}=0$. 
Furthermore if the metric also depends smoothly on time $t$ then
$$
 \sum_{j=0}^k \|D_t^j  P U\|_{H^r(\Omega)}\leq C_{r,k}\sum_{j=0}^k 
\|D_t^j   U\|_{H^r(\Omega)}.\tag 3.8
$$
This follows by induction in $k$ 
from commuting through time derivatives in (3.2):
$$
\triangle D_t^m q=-\sum_{j=0}^{m-1}{\text{$\binom{m}{j} $}}\pa_a\big( (D_t^{m-j} g^{ab})\pa_b D_t^j q\big)
+\pa_a \big( D_t^m U^a\big),\qquad D_t^m q\big|_{\pa\Omega}=0\tag 3.9
$$
which using (3.7) gives 
$\| D_t^m q\|_{H^{r+1}(\Omega)}\leq 
C_{r,m}\sum_{j=0}^{m-1}\| D_t^m q\|_{H^{r+1}(\Omega)}
+C_{r,m}\sum_{j=0}^m \| D_t^j U\|_{H^r(\Omega)}$.

For functions $f$ vanishing on the boundary
we define operators on divergence free vector fields
$$
A_f W^a=P \big(-g^{ab}\pa_b( (\pa_c f) W^c)\big), \tag 3.10
$$
$A_f$ is symmetric, i.e. $\langle U,A_f W\rangle=\langle A_f U,W\rangle$,
 since for $U$ and $W$ divergence free it follows from (3.3)
$$
\langle U,A_f W\rangle =- \int_\Omega U^a \pa_a\big( (\pa_c f) W^c\big)\,\kappa  dy
=\int_{\pa\Omega} (-\na_N f) U_N W_N \, \kappa dS,\qquad U_N= N_a U^a \tag 3.11
$$
If $p$ is the pressure in Euler's equations then normal operator $A$ in (2.26)
is 
$$
A=A_p\geq 0,\qquad\text{i.e.}\qquad \langle W,AW\rangle \geq 0,
\qquad\text{if}\qquad \na_N p\Big|_{{\pa\Omega}}\leq 0\tag 3.12
$$
which is true by our assumption (1.6). 
It follows from Cauchy Schwartz inequality that 
$$
|\langle U, A_{f\, p} W\rangle |\leq 
\langle U, A_{|f|p} U\rangle^{1/2}\langle W, A_{|f|p} W\rangle^{1/2}
\leq \|f\|_{L^\infty({\pa\Omega})} 
\langle U, A U\rangle^{1/2}\langle W, A W\rangle^{1/2}\tag 3.13
$$
since $\na_N (P)=f\na_N p$ on the boundary. The positivity properties 
(3.12) and (3.13)
are of fundamental importance to us. In particular, since $p$ vanishes on
the boundary so does $\dot{p}=D_t p$ and therefore 
$$
\dot{A}=A_{\dot{p}}\qquad\text{satisfies}\qquad 
|\langle W, \dot{A} W\rangle |
\leq \|\na_N \dot{p}/\na_N p\|_{L^\infty(\pa\Omega)}\langle W, A W\rangle 
\tag 3.14
$$
$\dot{A}$ is the time derivative of the operator $A$, considered 
as an operator with values in the one forms.

It follows from (3.10) and (3.5) that 
$
\| A_f W\|\leq \|\pa^2 f\|_{L^\infty(\Omega)} \|W\|+
\|\pa f\|_{L^\infty(\Omega)} \|\pa W\|.
$
However, $A_f$ acting on divergence free vector fields by (3.11) 
depends only on $\na_N f\big|_{\pa\Omega}$, i.e. $A_{\tilde{f}}=A_f$
if $\na_N \tilde{f}\big|_{\pa\Omega}=\na_N f\big|_{\pa\Omega}$. 
We can therefore replace $f$ by the Taylor expansion of order one in the distance 
to the boundary in polar coordinates multiplied by a smooth function that is one 
close to the boundary and vanishes close to the origin. It follows that 
$$
\|A_f W\|\leq C\sum_{S\in{\Cal S}}\|\na_N Sf\|_{L^\infty(\pa\Omega)}\|W\|
+ C\|\na_N f\|_{L^\infty(\pa\Omega)}(\|\pa W\|+\|W\|). \tag 3.15
$$
where ${\Cal S}$  is a set of vector fields that span the tangent space 
of the boundary, see section 6.

For two forms $\alpha$ we define bounded projected multiplication operators given by 
$$
M_\alpha {W}^a=P\big(g^{ab} \alpha_{bc} W^c\big),\qquad \quad
\|M_\alpha W\|\leq \|\alpha\|_{L^\infty(\Omega)}\|W\| .\qquad \tag 3.16
$$
In particular the operators in (2.27) and (2.28) are bounded projected multiplication
operators:  
$$
G=M_g ,\qquad C=M_\omega,\qquad \dot{G}=M_{\dot{g}}\tag 3.17
$$
where $g$ is the metric, $\omega$ the vorticity and $\dot{g}$ the 
time derivative of the metric.

\head 4. The lowest order energy estimate.\endhead 
Since $\det{(\pa x/\pa y)}=1$ it follows from introducing Lagrangian coordinates, 
that for a function $f$
$$
\int_\dt f \,dx=\int_\Omega f\, dy,\qquad \text{so}\quad 
\frac{d}{dt} \int_\dt f \, dx=\int_\dt D_t f \, dx\tag 4.1 
$$
We note that if $v$ is a solution if Euler's equations,
$D_t v_i=-\pa_i p$, and $p$ vanish on the boundary then 
$$
\frac{d}{dt} \int_\dt |V|^2 \, dx=2\int_\dt V^i D_t v_i\, dx
=-2\int_\dt V^i\pa_i p\, dx=2\int_\dt (\div V) p\, dx
-2\int_{\pdt} V_N p\, dS=0\tag 4.2
$$

We now want to obtain energy estimates for the linearized equations 
$$
L_1 W=\ddot{W}+ A W+\dot{G}\dot{W}-C\dot{W}=F\tag 4.3
$$
where $A$, $\dot{G}$ and $C$ are as in section 3 and 
$F$ is divergence free. Because of the unbounded but positive and symmetric 
operator $A$ there is an additional term in the energy:
$$
E=E(W)=\langle \dot{W},\dot{W}\rangle +\langle W, (A+I) W\rangle\tag 4.4
$$
where the inner product is given by (3.1).

Since $\langle\dot{W},\dot{W}\rangle = \int_\Omega g_{ab} \dot{W}^a\dot{W}^b\,
dy$ and $D_t \big(g_{ab}\dot{W}^a\dot{W}^b \big)=
\dot{g}_{ab}\dot{W}^a\dot{W}^b +2g_{ab} \dot{W}^a D_t \dot{W}^b$, we have 
$$
\frac{d}{dt} \langle \dot{W},\dot{W}\rangle 
=2\langle\dot{W}, D_t \dot{W}\rangle
+\langle\dot{W},\dot{G}\dot{W}\rangle \tag 4.5
$$
where $\dot{G}$ is given by (3.17). 
By (3.3) and (3.11) 
$\langle W, AW\rangle =-\int_{\Omega} W^a \pa_a\big( (\pa_c p) W^c\big)\, dy$,
and 
$$
D_t \big(W^a \pa_a \big((\pa_c p) W^c)\big)\big)
=\dot{W}^a\pa_a \big((\pa_c p) W^c)\big)
+W^a\pa_a \big((\pa_c p)\dot{W}^c\big)
+W^a\pa_a \big((\pa_c  D_t p) W^c\big).\tag 4.6
$$
Since $A$ is symmetric we get 
$$
\frac{d}{dt}\langle W, AW\rangle= 2\langle \dot{W},AW\rangle +
 \langle W,\dot{A}W\rangle\tag 4.7
$$
where $\dot{A} W^i=A_{\dot{p}} W^i$ is given by (3.10) with $f=\dot{p}=D_t p$. 
Hence 
$$\multline 
\frac{d}{dt} E(W)= 2\langle \dot{W}, \ddot{W}+AW+W\rangle 
+\langle\dot{W},\dot{G}\dot{W}\rangle +\langle W, \dot{A} W\rangle 
+\langle W,\dot{G} W\rangle \\
=2\langle \dot{W},L_1 W\rangle +2\langle \dot{W},W\rangle
-\langle \dot{W},\dot{G}\dot{W}\rangle+\langle W, \dot{A} W\rangle
+\langle W,\dot{G} W\rangle .
\endmultline \tag 4.8
$$
where we used that $\langle \dot{W},C\dot{W}\rangle$ vanishes since $C$ is antisymmetric. 
The operator $\dot{G}$ is bounded by (3.16)-(3.17) and 
$|\langle W,\dot{A} W\rangle|$ is bounded by (3.14) so
$$
|\dot{E}|\leq \Big(1+\|\dot{g}\|_{L^\infty(\Omega)}
+\|\na_N D_t p/\na_N p\|_{L^\infty(\pa\Omega)}\Big) E+2\sqrt{E} \|F\|. \tag 4.9
$$
With 
$n(t)=1+\|\dot{g}\|_{L^\infty(\Omega)}+\|\na_N D_t p/\na_N p\|_{L^\infty(\pa \Omega)}$
and $E_0=\sqrt{E}$ we hence have 
$$
E_0(t)\leq e^{\int_0^t n\, d\tau} 
\Big(E_0(0)+\int_0^t \|F\|\, d\tau\Big) . \tag 4.10
$$

\comment 
\head 5. Outline of the proof.\endhead
The strategy of the proof is to replace the unbounded normal operator by 
a sequence of bounded operators converging to it.
For the
equations with the bounded operators we have existence and we will show that we have 
uniform bounds for derivatives so the limit of the solutions will be regular.

First we will show {\it a priori} bounds assuming existence,
in such away that they 
generalize to give uniform bounds for the sequence of regularized equations. 
We want to commute tangential vector fields through (4.3) to get 
equations for higher order derivatives and use the energy estimate for these. 
The reason for only applying tangential vector fields,
constructed in section 7,  is that if $T$ is tangential 
then $p|_{\pa \Omega}\!=\!0$ implies that $Tp|_{\pa \Omega}\!=\!0$ 
and just as for the 
time derivative this will allow us to control the commutator with $A$. 
We will use Lie derivatives, see section 8, because they commute with 
the differential and the Lie derivative, 
with respect to a divergence free vector field, of a divergence free vector
 field is divergence free.
However, when applying these to the equation we will get 
a term with $AW$ or $D_t^2 W$ that we must bound. 
 Therefore we must also
control the energy of the higher time derivatives, so we must commute these 
through the equation as well. The initial conditions for the higher time
derivatives have to be obtained from the equation and it will be important to
turn these into an inhomogeneous term instead, see section 6.
In section 9 we calculate the commutators between the Lie derivatives with respect to
tangential vector fields and the normal operator and projection.
In section 10 we use this to commute the tangential vector fields through the
equation. This will give us control of
the energy of tangential components in section 11, see Lemma {10.}1. 
The additional components we
will get from that the curl satisfies a better equation since the 
curl of $A$ vanishes, see section 13, and that one can estimate all components 
from the tangential ones, the curl and the divergence, see section 11. 
The {\it a priori} estimates of all components are then given in Lemma {12.}2. 

Then we will, in section 14,
construct a bounded smoothed out normal operator $A^\varepsilon$ converging to the
normal operator $A$, as $\varepsilon\to 0$. It still has all the important 
properties like symmetry, positivity and almost exactly the same 
 commutators with tangential vector fields. Existence of solutions, $W^\varepsilon$,
for the smoothed out equation i.e. (4.3) with 
$A$ replaced by $A^\varepsilon$ follows directly from that $A^\varepsilon$ 
is bounded since its just an ordinary differential equation in Sobolev spaces. 
The basic energy estimate then gives us uniform bound for $\|W^\varepsilon\|$
so we can extract a weakly convergent subsequence converging to a limit $W$,
which then is a weak solution of the equation, see section 14. 
Therefore it only remains to show that that we have uniform bounds 
of the higher derivative of $W^\varepsilon$, see Lemma {15.}1,
so that the limit $W$ in in the 
Sobolev spaces $H^r$. 
The proof of this is almost the same 
as the proof of the {\it a priori} bounds, 
due to that the commutators with the 
smoothed out normal operator are almost the same as the ones with the 
normal operator. 
Then finally we have existence of solutions in Sobolev spaces,
 but because of the way we got them, by taking additional time derivatives,
we loose regularity compared to the initial conditions. 
However, once we have existence we can obtain better bounds for the solutions,
in Lemma {16.}1, 
and through an approximation procedure get the existence theorem stated above.
\endcomment

\head 5. Turning the initial conditions into an inhomogeneous 
divergence free term.\endhead 
As explained in the introduction we want to reduce the initial value
problem
$$
L_1 W=\ddot{W}+AW+\dot{G}\dot{W}-C\dot{W}=F,\qquad\quad W\big|_{t=0}=W_0,
\qquad \dot{W}\big|_{t=0}=W_1 \tag 5.1
$$
to the case of vanishing initial conditions and an inhomogeneous term $F$
that vanishes to any order as $t\to 0$.
This is achieved by subtracting off a power series solution in $t$ to (5.1):
$$
W^a_{0r}(t,y)=
\sum_{s=0}^{r+2} \frac{t^s}{s !}
W_s^a(y)\tag 5.2
$$
We note that if $W_s$ are divergence free it follows that $W_{0r}$ 
is divergence free. Here $W_0$ and $W_1$ are the initial conditions,
$W_2$ is obtained form the equation (5.1) at $t=0$:
$W_2=F-AW_0-\dot{G}W_1+C W_1$. Similarly, one gets higher order terms 
by first differentiating the equation with respect to time. 
It is clear that doing so we obtain an expression 
$D_t^{k+2} W=M_k(W,...,D_t^{k+1} W)+D_t^k F$
and from this we inductively define 
$W_{k+2}=M_k(W_0,...,W_{k+1})\big|_{t=0}+D_t^k F\big|_{t=0}$.
Here $M_k$ is some linear operator of order at most one and 
that is all we need to know. However, we are going to calculate
the explicit form of $M_k$ since we will do similar calculations later
on for other operators and this is a simple model case. 

Now it turns out that its easier to differentiate the corresponding operator 
with values in one forms; 
$$
\underline{L}_1 W_a=g_{ab} L_1 W^b
=g_{ab} \ddot{W}^b-\pa_a \big((\pa_c p) W^c\big)
+\pa_a q +(\dot{g}_{ab}-\omega_{ab}) \dot{W}^b=g_{ab} F^b\tag 5.3
$$
where $q$ is chosen so the last terms are divergence free,
and afterwards project the result to the divergence free vector fields. 
Let 
$$
q^s=D_t^s q\qquad 
p^s=D_t^{s} p,\qquad
g^s_{ab}=D_t^s g_{ab},\qquad
\omega^{s}_{ab}=D_t^{s} \omega_{ab}, \qquad F_s=D_t^s F\tag 5.4
$$
In general it follows from applying $D_t^r$ to (5.3), restricting 
to $t=0$ gives 
that 
$$
\sum_{s=0}^{r}{\tsize{\binom{r}{s}}\Big( g_{ab}^{r-s}
{W}_{s+2}^b-\pa_a \big((\pa_c p^{r-s} )W_s^c\big)\Big)
+\pa_a q^{r} }
+\sum_{s=0}^{r}{\tsize{\binom{r}{s}}}(g_{ab}^{{r-s}+1}-\omega_{ab}^{r-s} )
{W}_{s+1}^b=\sum_{s=0}^{r}{\tsize{\binom{r}{s}}}g_{ab}^{{r-s}}F_s\tag 5.5
$$
We now want to project each term onto divergence free vector fields. 
Let 
$$
A_s W^c=P \big(-g^{ab}\pa_b( (\pa_c p^s) W^c)\big),
\qquad G_s W^c=P\big( g^{ca} g^s_{ab} W^b\big), 
\qquad  C_s W^c=P\big( g^{ac}\omega^s_{ab} W^b\big)\tag 5.6
$$
We obtain
$$
W_{r+2}=
-\sum_{s=0}^{r-1}{\tsize{\binom{r}{s}} } G_{r-s} W_{s+2}
-\sum_{s=0}^{r}\tsize{\binom{r}{s}}
\big( G_{r-s+1} W_{s+1}-C_{r-s} W_{s+1}
+A_{r-s} W_s-G_{r-s}F_s \big)\tag 5.7
$$
This inductively defines $W_{r+2}$ from $W_0,..., W_{r+1}$. 
With $W_{0r}$ given by (5.2) we have hence achieved that 
$$
D_t^s\big( L_1 W_{0r}-F\big)\big|_{t=0}=0,\qquad\text{for}\quad s\leq r,
\qquad W_{0r}\big|_{t=0}=W_0,\quad \dot{W}_{0r}\big|_{t=0}=W_1\tag 5.8 
$$
Replacing $W$ by $W-W_{0r}$ and $F$ by $F-L_1 W_{0r}$ hence reduces (5.1) 
to the case of vanishing initial data and an inhomogeneous term that 
vanishes to any order $r$ as $t\to 0$. 

 We also note that if the initial data
are smooth  then we can construct a smooth approximate solution $\tilde{W}$
that satisfies the equation to all orders as $t\to 0$. This is obtained by multiplying
the $k^{th}$ term in (5.2) by a smooth cutoff $\chi(t/\varepsilon_k)$, to be chosen
below, and summing up the infinite series.
Here $\chi$ is smooth $\chi(s)=1$ for $|s|\leq 1/2$ and $\chi(s)=0$ for $|s|\geq 1$.
The sequence $\varepsilon_k>0$ can then be chosen small enough so that the series
converges in $C^m([0,T],H^m)$ for any $m$
if we take $(\|\tilde{W}_k\|_k+1)\varepsilon_k\leq 1/2$.

\head 6. Construction of the tangential vector fields.\endhead 
Let us now construct the tangential divergence free vector fields, that are 
time independent expressed in the Lagrangian coordinates, i.e. that commute 
with $D_t$:
$$
[D_t,T]=0.\tag 6.1 
$$
 This means that in the Lagrangian coordinates they are of the form
$T^a(y)\pa/\pa y^a$ and since  $\det{(\pa x/\pa y)}=1$ the divergence free condition 
is just
$$
\pa_a T^a=0.\tag 6.2 
$$
Since $\Omega$ is the unit ball in $\bold{R}^n$ 
the vector fields can be explicitly given. 
The vector fields 
$$
y^a\pa/\pa y^b-y^b\pa/\pa y^a\tag 6.3
$$ 
corresponding to rotations, span the tangent space of 
the boundary and are divergence free in the interior. 
Furthermore they span the tangent space of the level sets of the 
distance function from the boundary in the Lagrangian coordinates
$$
d(y)=\dist{(y,\pa \Omega)}=1-|y|\tag 6.4
$$ 
away from the origin $y\neq 0$. 
We will denote this set of vector fields by ${\Cal S}_0$ 
We also construct a set of divergence free 
vector fields that span the full tangent space 
at distance $d(y)\geq d_0$ and that are compactly supported 
in the interior at a fixed distance $d_0/2$ from the boundary. 
The basic one is 
$$
h(y^3,...,y^n)\Big(f(y^1)g^\prime(y^2)\pa/\pa y^1- 
f^\prime(y^1)g(y^2)\pa/\pa y^2\Big),\tag 6.5
$$ 
which is divergence free. 
Furthermore we can choose $f,g,h$ such that 
it is equal to $\pa/\pa y^1 $ when $|y^i|\leq 1/4$, for $i=1,...,n$
and so that it is $0$ when $|y^i|\geq 1/2$ for some $i$. 
In fact let $f$ and $g$ be smooth functions such that 
$f(s)=1$ when $|s|\leq 1/4$ and $f(s)=0$ when $|s|\geq 1/2$
and $g^\prime(s)=1$ when $|s|\leq 1/4$ and $g(s)=0$ when $|s|\geq 1/2$. 
Finally let $h(y^3,...,y^n)=f(y^3)\cdot\cdot\cdot f(y^n)$. 
By scaling, translation and rotation 
of these vector fields we can obviously construct a finite set 
of vector fields that span the tangent space when $d\geq d_0$ and are 
compactly supported in the set where $d\geq d_0/2$. 
We will denote this set of vector fields by ${\Cal S}_1$. 
Let ${\Cal S}={\Cal S}_0\cup {\Cal S}_1$ denote 
the family of tangential space vector fields  and let 
${\Cal T}={\Cal S}\cup \{D_t\}$ denote the family of space time 
tangential vector fields. 

Let the radial vector field be 
$$
R=c_1 y^a\pa/\pa y^a,\qquad c_1>0\tag 6.6
$$
Now, $\div R=n$ is not $0$ but for our purposes it suffices that 
it is constant since what we need is that if $\div W=0$ then
$\div {\Cal L}_R W=R\div W-W\div R=0$, where the Lie derivative ${\Cal L}_R$ is defined in the
next section. Let ${\Cal R}={\Cal S}\cup\{R\}$.
Note that ${\Cal R}$ span the full tangent space of the space everywhere. 
Let ${\Cal U}={\Cal S}\cup \{R\}\cup\{D_t\}$ denote the family of all the vector fields
construct above. 
Note also that the radial vector field commutes with the rotations;
$$
[R,S]=0,\qquad S\in {\Cal S}_0\tag 6.7
$$
Furthermore, the commutators of two vector fields in ${\Cal S}_0$ 
is just $\pm$ another vector field in ${\Cal S}_0$. 
Therefore, for $i=0,1$, let ${\Cal R}_i={\Cal S}_i\cup\{R\}$, ${\Cal T}_i={\Cal S}_i\cup\{D_t\}$
and ${\Cal U}_i={\Cal S}_i\cup\{R\}\cup\{ D_t\}$. 

Let ${\Cal U}=\{ U_i\}_{i=1}^M $ be some labeling of our family of vector fields. 
We will also use multindices $I=(i_1,...,i_r)$ of length $|I|=r$.
so $U^I=U_{i_1}\cdot\cdot\cdot U_{i_r}$ and
${\Cal L}_U^I={\Cal L}_{U_{i_1}}\cdot\cdot\cdot {\Cal L}_{U_{i_r}}$. 
Sometimes we will write ${\Cal L}_U^I$, where $U\in {\Cal S}_0$ or $I\in {\Cal S}_0$, 
meaning that $U_{i_k}\in {\Cal S}_0$ for all of the indices in $I$. 

Note also that the vector fields $U^a(y)\pa/\pa y^a$ expressed in the 
$x$ coordinates are given by $U^i\pa/\pa x^i$ where 
$U^i=U^a\pa x^i/\pa y^a$. We here use the convention that 
indices $a,....,f$ refers to the components in the Lagrangian frame and
indices $i,...,n$ refers to the components in the Eulerian frame.

\head 7. Lie derivatives.\endhead 
Let us now introduce the Lie derivative of the vector field $W$
with respect to the vector field $T$; 
$$
{\Cal L}_T W^a=TW^a -(\pa_c T^a) W^c\tag 7.1 
$$
We will only deal with Lie derivatives with respect to the vector fields $T$
constructed in the previous section.
For those vector fields $\div T=0$ so
$$
\div W=0\qquad \implies \qquad \div {\Cal L}_{T} W=T\div W-W\div T=0. \tag 7.2
$$
The Lie derivative of a one form is defined by
$$
{\Cal L}_T \alpha_a =T\alpha_a+(\pa_a T^c) \alpha_c.\tag 7.3
$$
The Lie derivatives also commute with exterior differentiation,
$[{\Cal L}_T, d]=0$ so if $q$ is a function, 
$$
{\Cal L}_T \pa_a q=\pa_a T q.\tag 7.4
$$
The Lie derivative of a two form is given by 
$$
{\Cal L}_T \beta_{ab} =T\beta_{ab}
+(\pa_a T^c) \beta_{cb}+(\pa_b T^c) \beta_{ac}.\tag 7.5
$$
Furthermore if $w$ is a one form and 
$\curl w_{ab}= dw_{ab}=\pa_a w_b-\pa_b w_a$ 
then since the Lie derivative commutes with exterior differentiation: 
$$
{\Cal L}_T \curl w_{ab}=\curl {\Cal L}_T w_{ab}.\tag 7.6
$$
We will also use that the Lie derivative satisfies 
Leibnitz rule, e.g. 
$$
{\Cal L}_T (\alpha_{c} W^c)=
({\Cal L}_T \alpha_{c}) W^c+\alpha_{c} {\Cal L}_T W^c,\qquad 
{\Cal L}_T (\beta_{ac} W^c)=
({\Cal L}_T\beta_{ac}) W^c+\beta_{ac} {\Cal L}_T W^c. \tag 7.7
$$
Furthermore, we will also treat $D_t$ as if it were a Lie derivative and we will set
$$
{\Cal L}_{D_t}=D_t.\tag 7.8
$$
Now of course this is not a space Lie derivative but rather could be interpreted 
as a space time Lie derivative in the domain $[0,T]\times \Omega$. 
But the important thing is that it satisfies all the properties of the 
other Lie derivatives we are considering, such as $\div W=0$
implies that $\div D_t W=0$ and $D_t \curl w=\curl D_t w$,
simply because it commutes with partial differentiation with respect to the
$y$ coordinates. The reason we use the notation (7.9) is that 
we will apply products of Lie derivatives and (7.9) and it is more 
efficient with the same notation.
Furthermore
$$
[{\Cal L}_{D_t},{\Cal L}_T]=0\tag 7.9
$$
this is because this quantity is ${\Cal L}_{[D_t,T]}$ and $[D_t,T]=0$ for the 
vector fields we are considering, or it follows from (7.1) and that 
$T^a=T^a(y)$ is independent of $t$.

\head 8. Commutators between Lie derivatives with respect to tangential vector
fields and the normal and multiplication operators.\endhead 
Note that the projection $P$ defined in section 3
almost commutes with the Lie derivative with respect to
tangential vector fields. 
In fact if denote the corresponding operator on
one forms by $\underline{P}$ 
$$
\underline{P}u_a=u_a-\pa_a q\tag 8.1 
$$
where $q$ is as in (3.2) and $u_a=g_{ab} U^b$, then 
$ {\Cal L}_{T} \underline{P} u_a = {\Cal L}_{T} u_a-\pa_a Tq$. Since 
$q=0$ on the boundary it follows that $Tq=0$ there so the last term vanishes if 
we project again: 
$$
\underline{P} ({\Cal L}_{T} \underline{P} u_a)=\underline{P}{\Cal L}_{T} u_a\tag 8.2
$$

We will need to calculate commutator between Lie derivatives with respect to 
tangential
vector fields $T$ and the operator $A_f$ defined in section 3.
Let $\underline{A}_f$ denote the 
corresponding operator taking a vector field to the one form 
$$
\underline{A}_f W_a=g_{ab} A_f W^b=-\pa_a\big( (\pa_c f) W^c-q\big),\tag 8.3
$$
Then since 
$$
{\Cal L}_T \pa_a\big( (\pa_c f) W^c\big)=
 \pa_a\big( (\pa_c T f) W^k\big)+ \pa_a\big( (\pa_c f) {\Cal L}_T  W^c\big)\tag 8.4
$$
it follows from (8.2) 
$$
\underline{P} {\Cal L}_{T} \underline{A}_f W_a  = \underline{A}_f {\Cal L}_T W_a
+\underline{A}_{Tf} W_a\tag 8.5
$$
Note that if $f=p$ then it follows from (3.13) that the commutator is 
lower order. In fact $p=0$ on the boundary implies that $Tp=0$ on the boundary 
if $T$ is a tangential vector field. 
Since $\na_N p\neq 0$ it follows $Tp/p$ is a continuous function 
that is equal to $\na_N Tp/\na_N p$ on the boundary. Hence by (3.13) 
$$
|\langle W, A_{Tp} W\rangle |\leq \| \na_N Tp/\na_N p\|_{L^\infty({\pa\Omega})} 
\langle W, A W\rangle\tag 8.6
$$

In view of (8.2) it follows that the multiplication operator $M_\alpha$, 
defined by (3.16) in section 3,  
satisfies the commutator relation
$$
\underline{P}{\Cal L}_T \underline{M}_{\,\alpha} W=\underline{M}_{\,\alpha} {\Cal L}_T W
+\underline{M}_{\,{\Cal L}_T\alpha} W,\qquad\text{where}\qquad 
\underline{M}_\alpha W_a =g_{ab} M_\alpha W^b \tag 8.7
$$
for a two form $\alpha$. Let 
$$
G_T=M_{g^T},\qquad  g^T_{ab}={\Cal L}_T g_{ab},\qquad\quad
C_T=M_{\omega^T}, \omega^T={\Cal L}_T \omega \tag 8.8
$$
We will also use special notation for the time derivatives of $G$:
$$
\dot{G} =G_{D_t}=M_{\dot{g}},\qquad \dot{g}_{ab}={D_t}g_{ab} \tag 8.9
$$
and of $A$
$$
A_T=A_{Tp},\qquad \dot{A}=A_{D_t}=A_{D_t p}\tag 8.10
$$
In the following sections 
we will commute through products of vector fields 
${\Cal L}_T^I={\Cal L}_{T_{i_1}}\cdot\cdot\cdot {\Cal L}_{T_{i_r}}$
where $I=(i_1,...,i_r)$ and we will use the notation  
$$
A_I=A_{T^I p},\qquad G_I=M_{g^I},\qquad C_I=M_{\omega^I}\qquad
\qquad \dot{G}_I=M_{\dot{g}^I}\tag 8.11 
$$
where $g^I_{ab}={\Cal L}_T^I g_{ab}$,  
$\dot{g}^I_{ab}={\Cal L}_T^I {D_t} g_{ab}$
and $\omega^I_{ab}={\Cal L}_T^I \omega_{ab}$.

\head 9. Commutators between the linearized equation and Lie derivatives with
respect to tangential vector fields.\endhead 
We are now ready to commute tangential vector fields through the linearized equation 
and in the next section get the higher order energy estimates of tangential derivatives.
Let $T\in{\Cal T}$ be a  tangential vector fields and recall 
that $[{\Cal L}_T,{D_t}]=0$ and that 
if $W$ are divergence free then so is ${\Cal L}_T W$.
Let us now apply Lie derivatives
 ${\Cal L}_T^I={\Cal L}_{T_{i_1}}\cdot \cdot \cdot{\Cal L}_{T_{i_r}}$, 
where $I=(i_1,...,i_r)$ is a multi index, to the linearized equation (2.19) with an inhomogeneous 
divergence free term $F$ vanishing to order $r$ as $t\to 0$: 
$$
g_{ab}   \ddot{W}^b-\pa_a \big((\pa_c p) W^c\big)
=-\pa_a q -(\dot{g}_{ca}-\omega_{ca}) \dot{W}^c+g_{ab} F^b,
\qquad W\big|_{t=0}=\dot{W}\big|_{t=0}=0.\tag 9.1
$$
which yields 
$$\multline 
c_{I_1 I_2}^{\,I}( {\Cal L}_T^{I_1} g_{ab} )
  {\Cal L}_T^{I_2}\ddot{W}^b-
c_{I_1 I_2}^{\,I}\pa_a \big((\pa_c T^{I_1} p)  {\Cal L}_T^{I_2} W^c\big)\\
=-\pa_a T^{I} q -2 c_{I_1 I_2}^{\,I} \big( {\Cal L}_T^{I_1}(\dot{g}_{ca}-\omega_{ca}) \big) 
{\Cal L}_T^{I_2}\dot{W}^c+c_{I_1 I_2}^{\,I}
( {\Cal L}_T^{I_1} g_{ab} ){\Cal L}_T^{I_2} F^b
\endmultline 
\tag 9.2
$$
where we sum over all $I_1+I_2=I$ and $c_{I_1 I_2}^{\,I}=1$. 
Let us introduce some new notation
$$
W_I={\Cal L}_T^I W,\quad {F}_I={\Cal L}_T^I {F}\quad
g^I_{ab}={\Cal L}_T^I g_{ab}
\quad
\omega^I_{ab}={\Cal L}_T^I\omega_{ab} ,
\quad 
p_I=T^I p,\quad q_I=T^I q \tag 9.3 
$$
and $\dot{g}^I_{ab}= {D_t}{\Cal L}_T^I g_{ab}$,
$\dot{W}_I={D_t} W_I$ etc. With this notation (9.2) becomes
$$
c_{I_1 I_2}^{\,I} g^{I_1}_{ab} \ddot{W}^b_{I_2} 
-c_{I_1 I_2}^{\,I} \pa_a\big( ( \pa_c p_{I_1}) W_{I_2}^c\big)=-\pa_a q_{I} 
-c_{I_1 I_2}^{\,I} \big( \dot{g}_{ab}^{I_1} -\omega^{I_1}_{ab}\big) \dot{W}^b_{I_2}
+c_{I_1 I_2}^{\,I} g^{I_1}_{ab} F_{I_2}^b\tag 9.4 
$$
Let us now project each term onto divergence free vector fields and
also introduce some notation for the resulting operators
$$
A_I W^a =A_{T^I p} W^a ,\qquad G_I W^a=P\big( g^{ac} g_{cb}^I W^b)\tag 9.5
$$
and 
$$
\dot{G}_I W^a=P\big( g^{ac}\dot{g}_{cb}^I W^b),\qquad 
C_I W^a= P( g^{ac}\omega_{cb}^I W^b)\tag 9.6
$$
From now on we set $\tilde{c}^{\, I_1 I_2}_I=c_{I_1 I_2}^{\,I}$ when $I_2\neq I$
and $\tilde{c}^{\, I_1 I_2}_I=0$ if $I_2=I$. 
Projecting each term onto divergence free vector fields we can now write (9.4) as
$$
L_1 W_I=\ddot{W}_I+ A W_I+\dot{G}\dot{W}_I-C \dot{W}_I=F_I
-\tilde{c}^{\,I_1 I_2}_I\big( A_{{I_1}} W_{I_2} +\dot{G}_{I_1}
\dot{W}_{I_2}- C_{I_1} \dot{W}_{I_2}
+G_{I_1} \ddot{W}_{I_2}+G_{I_1} F_{I_2}\big)
\tag 9.7
$$
Here $G_{J}$, $\dot{G}_J$ and $C_J$ are all bounded operators. By (3.16)-(3.17):
$$
\|G_J W\|\leq \|{\Cal L}_T^J  g\|_{L^\infty(\Omega)} \|W\|,
\qquad \|C_J W\|\leq \|{\Cal L}_T^J   \omega\|_{L^\infty(\Omega)} \|W\|\tag 9.8 
$$
The terms $G_{I_1} \ddot{W}_{I_2}$ 
are easy to take care of by also including time derivatives
up to highest order in our estimates 
since $|I_2|\leq |I|-1$. 
$A W_I$ itself will be included in the higher order 
energy, which is just going to be a sum of terms of the form (4.4) 
with $W$ replaced by $W_I$ for $|I|\leq r$.
However, we also have to deal with $A_{I_1}W_{I_2} $ since $A_{I_1}$ is an operator 
of order $1$. Since $|I_2|\leq |I|-1\leq r-1$ in the terms
$A_{{I_1}} W_{I_2}$ and since the energy will give us $\dot{W}_I$
for all $|I|\leq r$ we in particular will have an estimate for 
$\ddot{W}_{I_2}$ which, using the equation (9.7),
up to terms of lower order is $-A W_{I_2}$.
Since $A_{J}=A_{T^J p}$ it follows from
(3.13) that 
$$
|\langle U,A_{J} W\rangle |\leq \| \na_N T^J p/ \na_N p \|_{L^\infty(\pa\Omega)}
\langle U, AU\rangle^{1/2}\langle W, AW\rangle^{1/2},\tag 9.9 
$$
However this does not imply that the norm of $A_J$ is bounded by the norm of $A$.
Therefore we have to deal with these terms with $A_{I_1}$ in an indirect way,
by including them in the energy and using (9.9). 

\head 10. The {\it a priori} energy bounds for tangential derivatives.\endhead
To obtain estimates for higher derivatives we apply tangential vector fields
to the equation and get similar equations for higher derivatives.
However, there are a some commutators coming up that we have to deal with.
One can be dealt with by
adding a lower order term to the energy and another commutator one deals
with by also considering higher time derivatives.
The main point is however that commutators with the normal operator can
be controlled by the normal operator through (9.9).
Let $W_T={\Cal L}_T W$, $F_T={\Cal L}_T F$ and let 
$G_T$, $C_T$ and $A_T$ be as in (8.8) and (8.10). By (9.7)
$$
L_1 W_T=F_T-A_T W-\dot{G}_T W+C_T \dot{W}
+G_T \ddot{W}+G_T F\tag 10.1
$$
The terms one has to deal with are $A_T W$ and $G_T \ddot{W}$.
Let $E=E(W)$, where $E(W)$ is given by (4.4),  
$$
E_T=E(W_T)=\langle \dot{W}_T,\dot{W}_T\rangle +
\langle W_T, (A+I)W_T\rangle,\qquad\text{and}\qquad 
D_T=2\langle W_T,A_T W\rangle.\tag 10.2
$$
$D_T$ is lower order compared to $E_T$ since by (9.9) 
it is bounded by a constant times $\sqrt{E_T}\sqrt{E}$ and we already have an
estimate for $E$ in (4.10). We will add $D_T$ to the energy $E_T$ to
pick up the commutator $A_T$ between ${\Cal L}_T$ and $A$.
By (4.8)
$$\multline
\dot{E}_T+\dot{D}_T=2\langle \dot{W}_T,L_1 W_T\rangle
+2\langle \dot{W}_T,W_T\rangle
-\langle \dot{W}_T,\dot{G}\dot{W}_T\rangle+\langle W_T, \dot{A} W_T\rangle
+\langle W_T,\dot{G} W_T\rangle
\\
+2\langle \dot{W}_T,A_T {W}\rangle 
+2\langle W_T,A_T\dot{W}\rangle 
+2\langle W_T,\dot{A}_T W\rangle\\
=2\langle W_T,A_T\dot{W}\rangle+2\langle \dot{W}_T ,G_T\ddot{W}\rangle
+2\langle \dot{W}_T, F+G_T F\rangle\\
+ 2\langle \dot{W}_T,-\dot{G}_T W+C_T \dot{W}+W_T\rangle 
-\langle \dot{W}_T,\dot{G}\dot{W}_T\rangle+\langle W_T, \dot{A} W_T\rangle
+\langle W_T,\dot{G} W_T\rangle
+2\langle W_T,\dot{A}_T W\rangle
\endmultline\tag 10.3
$$
Here, the terms on the last row are bounded by $E_T$ and $E$ using (9.8) and (9.9).
The only terms that remains to control are
$2\langle \dot{W}_T,G_T \ddot{W}\rangle$ and $2\langle W_T,A_T\dot{W}\rangle$.
These terms are controlled by simultaneously consider one more time
derivative, i.e. if $T=D_t$, and estimate energies for these.

Let us now define higher order energies. Let 
$$
E_I=E(W_I)= \langle \dot{W}_I,\dot{W}_I\rangle +
\langle W_I, (A+I)W_I\rangle,\qquad W_I={\Cal L}_T^I W.\tag 10.4
$$
With notation as in the previous section we have by (4.8) and (9.7) 
$$\multline 
\dot{E}_I=2\langle \dot{W}_I,{D_t} \dot{W}_I+ A
W_I+\dot{G}\dot{W}_I-C \dot{W}_I \rangle\\
 +2\langle \dot{W}_I,W_I\rangle
-\langle \dot{W}_I,\dot{G}\dot{W}_I\rangle+\langle W_I, \dot{A} W_I\rangle
+\langle W_I,\dot{G} W_I\rangle \\
=-2\tilde{c}^{\,I_1 I_2}_I\Big(\langle \dot{W}_I,A_{{I_1}} W_{I_2}\rangle
+\langle \dot{W}_I,\dot{G}_{I_1} \dot{W}_{I_2}\rangle - \langle \dot{W}_I,C_{I_1}  \dot{W}_{I_2}\rangle
+\langle \dot{W}_I,G_{I_1} \ddot{W}_{I_2}\rangle+\langle \dot{W}_I,G_{I_1} F_{I_2}
\rangle \Big)\\
+2\langle \dot{W}_I , F_I\rangle +2\langle \dot{W}_I,W_I\rangle
-\langle \dot{W}_I,\dot{G}\dot{W}_I\rangle+\langle W_I, \dot{A} W_I\rangle 
+\langle W_I,\dot{G} W_I\rangle 
\endmultline \tag 10.5
$$
To deal with the term $ \langle \dot{W}_I,A_{I_1} W_{I_2}\rangle $ 
we introduce
$$
D_I=2\, \tilde{c}^{\,I_1 I_2}_I\langle W_I , A_{I_1} W_{I_2}\rangle \tag 10.6
$$
Then 
$$
\dot{D}_I=2\, \tilde{c}^{\,I_1 I_2}_I\Big(\langle \dot{W}_I , A_{I_1} W_{I_2}\rangle
+\langle {W}_I , A_{I_1} \dot{W}_{I_2}\rangle
+\langle {W}_I , \dot{A}_{I_1} W_{I_2}\rangle\Big)\tag 10.7
$$
and hence 
$$\multline
\dot{E}_I+\dot{D}_I=\\
-2\, \tilde{c}^{\,I_1 I_2}_I\Big(-\langle {W}_I , A_{I_1} \dot{W}_{I_2}\rangle
-\langle {W}_I , \dot{A}_{I_1} W_{I_2}\rangle
+\langle \dot{W}_I,\dot{G}_{I_1} \dot{W}_{I_2}\rangle - \langle \dot{W}_I,C_{I_1}  \dot{W}_{I_2}\rangle
+\langle \dot{W}_I,G_{I_1} \ddot{W}_{I_2}\rangle+\langle \dot{W}_I,G_{I_1} F_{I_2}
\rangle \Big)\\
+2\langle \dot{W}_I , F_I\rangle +2\langle \dot{W}_I,W_I\rangle
-\langle \dot{W}_I,\dot{G}\dot{W}_I\rangle+\langle W_I, \dot{A} W_I\rangle 
+\langle W_I,\dot{G} W_I\rangle 
\endmultline \tag 10.8
$$
We have hence replaced the bad term by two terms that we can control by (9.9).
Furthermore, we can also bound $D_I$ itself using (9.9). 

For a two form $\alpha$ and a function $q$ vanishing on the boundary let 
$$ 
\|\alpha\|_\infty=\|\, |\alpha|\, \|_{L^\infty(\Omega)},
\qquad
\|\pa q\|_{\infty,\, p^{-1}}=\|\na_N q/\na_N p\|_{L^\infty(\pa \Omega)}
\leq \|\pa q\|_\infty/c_0,\tag 10.9 
$$
and for a vector fields $W$
let 
$$
\langle W\rangle_A=\langle W,A W\rangle^{1/2},\qquad 
\|W\|=\langle W,W\rangle^{1/2}.\tag 10.10
$$

With this notation it now follows from (10.8) and (9.8)-(9.9) that 
$$\multline 
\dot{E}_I+\dot{D}_I\leq 2\langle W_I\rangle_A
 \tilde{c}^{\,I_1 I_2}_I \Big(  \|\pa {p}_{I_1}\|_{\infty,\, p^{-1}} \langle \dot{W}_{I_2}\rangle_A
+\|\pa \dot{p}_{I_1}\|_{\infty, \, p^{-1}}\langle W_{I_2}\rangle_A\Big)\\
+2\|\dot{W}_I\|\, \tilde{c}^{\,I_1 I_2}_I\Big(  (\|\dot{{g}}^{I_1}\|_\infty+\|\omega^{I_1}\|_\infty)
 \|\dot{W}_{I_2}\|+ \|{{g}}^{I_1}\|_\infty \|\ddot{W}_{I_2}\|
+ \|{{g}}^{I_1}\|_\infty \|F_{I_2}\|\Big)\\
+\|\dot{W}_I\|\big( 2\|F_I\|+2\|W_I\|+\|\dot{{g}}\|_\infty \|\dot{W}_I\|\big)
+\|W_I\|\, \|\dot{{g}}\|_{\infty}\|{W}_I\|+
\langle W_I\rangle_A \|\pa \dot{p}\|_{\infty,\, p^{-1}} \langle W_I\rangle_A
\endmultline \tag 10.11
$$
Furthermore 
$$
|D_I|\leq 2\langle W_I\rangle_A \, \tilde{c}^{\,I_1 I_2}_I 
\| \pa p_{I_1}\|_{\infty,\, p^{-1}} \langle W_{I_2}\rangle_A \tag 10.12
$$

\demo{Definition {10.}1} For ${\Cal V}$ any of our families of vector fields let 
$$
E_{s}^{\Cal V}= \sum_{|I|\leq s,\,I\in{\Cal V}} \sqrt{E_I},\qquad \quad 
\|W\|_{s}^{\Cal V}=\sum_{|I|\leq s,\,I\in{\Cal V}} \|{\Cal L}_{T}^I W\|\tag 10.13
$$
where $E_I$ is given by (10.4). 
For a two form $\alpha$ and a function $q$ vanishing on the boundary let 
$$
\|\alpha \|_{s,\infty}^{\Cal V}=\sum_{|J|\leq s,\,
J\in {\Cal V}}\|{\Cal L}_T^J \alpha \|_\infty,
\qquad \|\pa q\|_{s, \infty,\, p^{-1}}^{\Cal V}=\sum_{|J|\leq s, \, J\in {\Cal V}}
\|\pa T^J q \|_{\infty,\, p^{-1}}
\tag 10.14
$$
where the norms are given by (10.9). Furthermore, let 
$$
n_s^{\Cal V}=\|\dot{g}\|_{s,\infty}^{\Cal V}+\|\omega\|_{s,\infty}^{\Cal V}
+\|{g}\|_{s+1,\infty}^{\Cal V}+\|\pa p\|_{s+1,\infty,p^{-1}}^{\Cal V} 
+\|\pa \dot{p}\|_{s,\infty,p^{-1}}^{\Cal V} \tag 10.15
$$
\enddemo
If $I\in{\Cal T}$ and 
 $|I|=r$ then with the notation in Definition {10.}1 we obtain from
(10.11) and (10.12): 
$$
|\dot{E}_I+\dot{D}_I|\leq C E_r^{\Cal T}
\sum_{s=0}^r n_s^{\Cal T} \big(E_{r-s}^{\Cal T}+\|F\|_{r-s}^{\Cal T}\big),
\qquad\qquad |D_I|\leq CE_r^{\Cal T} \sum_{s=0}^{r-1} n_s^{\Cal T} E_{r-1-s}^{\Cal T}
\tag 10.16
$$
If we integrate the first inequality from $0$ to $t$ using that $E_I(0)=D_I(0)=0$
and the second inequality we get with a constant depending on 
$\overline{n}_r=\sup_{0\leq \tau\leq T} n_{s}^{\Cal T}(\tau)$ 
$$
E_I\leq C E_r^{\Cal T} E_{r-1}^{\Cal T} 
+C\int_0^t E_{r}^{\Cal T} \big(  E_{r}^{\Cal T} +\|F\|_r^{\Cal T}\big)\, d\tau
\tag 10.17 
$$
If we sum over $|I|\leq r$ and divide by 
$ \overline{E}_r(t)=\sup_{0\leq \tau \leq t}
 E_r^{\Cal T}(\tau)$
we get for some other constant 
$$
\overline{E}_r\leq C\overline{E}_{r-1}
+C\int_0^t\big( \overline{E}_r+\|F\|_{r}^{\Cal T}\big)\, d\tau.\tag 10.18 
$$
Hence with $M_r(t)=\int_0^{t} \overline{E}_r\, d\tau$, we get 
$$
\frac{d{M}_r}{dt}-CM_r \leq C \overline{E}_{r-1} + C\int_0^t \|F\|_r^{\Cal T} \, d\tau \tag 10.19 
$$
Multiplying by the integrating factor $e^{-Ct}$ and integrating from $0$ to $t$ 
we see that $M_r$ is bounded by some constant depending on $t\leq T$ 
times the right hand side and hence it follows that for some other constant
$$
\overline{E}_r\leq C\overline{E}_{r-1} +C\int_0^t\|F\|_r^{\Cal T}\, d\tau\tag 10.20
$$
Since we already proved a bound for $\overline{E}_0$ in (4.10) it inductively follows
that: 

\proclaim{Lemma {10.}1} Suppose that $x, p\in C^{r+2}([0,T]\times\Omega)$, 
$p\big|_{\pa\Omega}=0$,
$\na_N p\big|_{\pa\Omega}\leq -c_0<0$ and $\div V=0$, where $V=D_t x$.
Suppose that $W$ is a solution of (9.1) where $F$ is divergence free and 
vanishing to order $r$ as $t\to 0$. Let $E^{\Cal T}_s$ be defined by (10.14). 
Then there is a constant $C$ depending only on the norm 
of $(x,p)$, a lower bound for $c_0$ and an upper bound for $T$, 
 such that if $E^{\Cal T}_s(0)=0$, for $s\leq r$, then 
$$
  E^{\Cal T}_r(t)\leq C
\int_0^t   \|F\|_{r}^{\Cal T}\, d\tau ,\qquad\text{for}\quad 0\leq t\leq T
\tag 10.21
$$
\endproclaim

\head 11. Estimates of derivatives of a vector field in terms of the 
curl, the divergence and tangential derivatives.\endhead 
In this section we show that derivatives of vector fields can be estimated
by derivatives of the curl, the divergence and tangential derivatives. 
First we prove the basic estimate in the Euclidean coordinates in Lemma {11.}1 below.
This estimate it is not invariant and so in Lemma {11.}2 we express it in terms
of Lie derivatives which is invariant. 
\proclaim{Lemma {11.}1} We have 
$$
|\pa {\alpha}|\leq 
C_n\big( |\curl {\alpha}|+|\div {\alpha}|+\tsize{\sum_{S\in{\Cal S}}} |S{\alpha}|\big),
\qquad \curl \alpha_{ij}=\pa_i \alpha_j-\pa_j \alpha_i\qquad
\div \alpha=\delta^{ij}\pa_i \alpha_j \tag 11.1
$$
for a one form $\alpha_i$ in the Eulerian frame, 
where $C_n$ only depends on the dimension $n$. 
Here the norms are the Euclidean norms,
$|\pa\alpha|=\sqrt{\sum_{i,j=1}^n |\pa_i\alpha_j|^2}$. 
\endproclaim
\demo{Proof of Lemma {11.}1. }
Since ${\Cal S}$ span the full tangent space 
in the interior when the distance to the boundary 
$d(y)\geq d_0$ we may assume that $d(y)<d_0$.
Let $\Omega^a=\{y;\, d(y)>a\}$ and let $\dt^a$ be the image of this set
under mapping $y\to x(t,y)$. 
Let ${N}$ the exterior unit normal to ${\pdt}^a$. 
Then $q^{ij}=\delta^{ij}-{N}^i {N}^j$ is the inverse of the tangential
metric. Since the tangential vector fields 
span the tangent space of the level sets of the distance function we have 
 $q^{ij} a_i a_j\leq C\sum_{S\in {\Cal S}} S^i S^j a_i a_j$,
where here $S^i=S^a\pa x^i/\pa y^a$.  
We claim that for any two tensor $\beta_{ij}$:
$$
\delta^{ij}\delta^{kl} \beta_{ki}\beta_{lj}\leq 
C_n \big(\delta^{ij}q^{kl}\beta_{ki}\beta_{lj}+|\hat{\beta}|^2+(\tr{\beta})^2\big)
\tag 11.2
$$
where $\hat{\beta}_{ij}=\beta_{ij}-\beta_{ji}$ is the antisymmetric part and 
$\tr \beta=\delta^{ij} \beta_{ij}$ is the trace. To prove (11.2) we may assume
 that $\beta$ is symmetric and
traceless. Writing $\delta^{ij}=q^{ij}+ {N}^i  {N}^j$ we see that the estimate for such
tensors follows from the estimate 
$ {N}^i  {N}^j  {N}^k  {N}^l \beta_{ki}
\beta_{lj}
=( {N}^i  {N}^k \beta_{ki})^2=(q^{ik} \beta_{ki})^2\leq 
n q^{ij} q^{kl} \beta_{ki} \beta_{lj}$. (This inequality just says that
$(\tr(Q\beta))^2 \leq n \tr(Q\beta Q\beta)$ which is obvious if 
one writes it out and use the
symmetry. ) 
\qed\enddemo 

The inequality (11.1) is not invariant under changes of coordinates so we 
want to replace it by an inequality that is, so we can get an inequality
that holds also in the Lagrangian frame. 
After that we want to derive higher order versions of it as well.
The divergence and the curl are invariant but the other terms are not.
There are two ways to make these terms invariant.  
One is to replace the differentiation by covariant differentiation and the 
other is to replace it by Lie derivatives with respect to the 
our family of vector fields in section 6. 
Both ways will result in a lower order term just involving 
the norm of the one form itself multiplied by a constant which 
depends on two derivatives of the coordinates. 
\comment
In the inequalities below, $C$ will stand for a constant that depends only 
on bounds for a finite number of derivatives of the vector fields in ${\Cal U}$
expressed in the Lagrangian coordinates $\pa^\alpha U^a(y)/\pa y^\alpha$. 
\endcomment 

\demo{Definition {11.}1} Let $c_1$ be a constant such that
$$
\sum_{a,b}\big(|g_{ab}|+|g^{ab}|\big)\leq c_1^2,
\qquad\quad |\pa x/\pa y|^2+|\pa y/\pa x|^2\leq c_1^2\tag 11.3
$$
and let $K_1$ denote a continuous function of $c_1$.
\enddemo
We note that the bound for the Jacobian of the coordinate and its inverse 
follows from the bound for the metric and its inverse and the bound 
for the Jacobian and its inverse implies an equivalent bound for the metric
and its inverse with $c_1^2$ multiplied by $n$. All our constants in what follows
in this section will depend on a bound for $c_1$ and we will denote such a 
constants by $K_1$.
\proclaim{Lemma {11.}2} In the Lagrangian frame we have, with 
$\underline{W}_a=g_{ab} W^b$, 
$$\align 
|{\Cal L}_U W|&\leq K_1\Big(|\curl\,\underline{W}\,|+|\div {W}|
+\tsize{\sum_{S\in{\Cal S}}} |{\Cal L}_S W|+[g]_1|W|\Big), 
\qquad U\in{\Cal R},\tag 11.4\\
|{\Cal L}_U W|&\leq K_1\Big(|\curl\, \underline{W}\,|+|\div {W}|
+\tsize{\sum_{T\in{\Cal T}}} |{\Cal L}_T W|+[g]_1|W|\Big), 
\qquad U\in{\Cal U},\tag 11.5
\endalign 
$$
where $[g]_1=1+|\pa g|$. 
Furthermore 
$$
|\pa W|\leq K_1\Big(|{\Cal L}_R W|
+\tsize{\sum_{S\in{\Cal S}}} |{\Cal L}_S W|+ |W|\Big)
\tag 11.6
$$
When $d(y)\leq d_0$ we may replace the sums over ${\Cal S}$
by the sums over ${\Cal S}_0$ and 
the sum over ${\Cal T}$
by the sum over ${\Cal T}_0$. 
\endproclaim
\demo{Proof of Lemma {11.2} } (11.5) follows 
directly from (11.4) by adding the time derivative to the right hand side.
We will show that (11.4) in the Eulerian frame follows from (11.1)
and then it follows directly
that (11.4) holds in Lagrangian frame as well since everything is invariant. 
Let $Z^i=\delta^{ij}\alpha_j$. Then
${\Cal L}_U Z^i=U Z^i-(\pa_k \tilde{U}^i) Z^k$,
where $\tilde{U}^i=U^a \pa x^i/\pa y^a$ are the components of the vector field
$U$ expressed in the Eulerian frame. Now transforming to the Lagrangian frame,
partial differentiation becomes covariant differentiation.
$ (\pa_k \tilde{U}^i)(\pa x^k/\pa y^a)(\pa y^b/\pa x^i) =\na_a U^b$,
where $\na_a U^b=\pa_a U^b+\Gamma_{a\,\, \, c}^{\,\,\, \,  b} U^c$,
and $\Gamma_{ab}^{\,\,\,\,\,\, c}={g^{cd}}
\left(\pa_a{g_{bd}}+\pa_b{ g_{ad}}-\pa_d{g_{ab}}\right)/2
=({ \partial y^c}/{\partial x^i}){\pa_a\pa_b x^i}$ are the Christoffel symbols.
Since $|\pa_a U^b|\leq C$ it follows that 
$|\pa_k \tilde{U}^i|\leq C[g]_1$.
That we may replace ${\Cal S}$ by ${\Cal S}_0$ close to the boundary 
follows from the proof of Lemma {11.}1. (11.6) follows since 
${\Cal R}$ span the tangent space and
$|{\Cal L }_{U} W^a-UW^a|=|(\pa_c U^a)W^c|\leq C|W|$. 
\qed\enddemo

We are now going to derive higher order versions of the inequality in Lemma {11.}2.
We want to apply the lemma to $W$ replaced by ${\Cal L}_U^J W $. 
Then in our applications the divergence term vanishes 
and as we shall see later on we will be 
able to control the curl of $({\Cal L}_U^J \underline{W})_a = {\Cal L}_U^J( g_{ab} W^b)$ 
which however is not the same as 
the curl of $(\underline{ {\Cal L}_U^J W})_a=g_{ab} {\Cal L}_U^J W^b $ but the difference is lower order and can
be easily estimated. Let us first introduce some notation:
\demo{Definition {11.}2} Let $\beta$ be a function, a one or two form or
 vector field, 
let ${\Cal V}$ be any of our families of vector fields and set
$$
|\beta|_{s}^{{\Cal V}}=\!\!\!\!\!\!\!\!
\sum_{|J|\leq s,\, J\in {\Cal V}}\!\!\!\!\!\! |{\Cal L}_S^J \beta|,
\qquad\qquad
[\beta]_{u}^{\Cal V} =\!\!\!\!\!\!\!\!\sum_{s_1+...+s_{k}\leq u,\, s_i\geq 1}\!\!\!\!\!\!\!\!
 |\beta|_{s_1}^{\Cal V}\!\!\cdot\cdot\cdot |\beta|_{s_{k}}^{\Cal V},
\qquad\qquad [\beta]_{0}^{\Cal V} =1. \tag 11.7
$$
\enddemo
In particular $|\beta|_{r}^{\Cal R}$ is equivalent to 
$\sum_{|\alpha|\leq r}|\pa^\alpha_y \beta|$ and 
$|\beta|_{r}^{\Cal U}$ is equivalent to 
$\sum_{|\alpha|+k\leq r}|D_t^k \pa^\alpha_y \beta|$. 
\proclaim{Lemma {11.}3} With the convention that 
$|\curl \underline{W}|_{-1}^{\Cal V}=|\div W|_{-1}^{\Cal V}=0$ we have 
$$\align 
|W|_{r}^{{\Cal R}}&\leq 
 K_1 \big(\,|\curl \underline{W}|_{r-1}^{\Cal R}+|\div W|_{r-1}^{\Cal R}
+ |W|_{r}^{\Cal S}
+\sum_{s=1}^{r}|g|_{s}^{\Cal R} |W|_{r-s}^{\Cal R}\big),
\tag 11.8\\
|W|_{r}^{{\Cal R}}&\leq 
 K_1 \sum_{s=0}^{r}\,\,
[g]^{\,\Cal R}_{s}\big(\,|\curl \underline{W}|_{r-1-s}^{\Cal R}+|\div W|_{r-1-s}^{\Cal R}
+|W|_{r-s}^{\Cal S}\big). \rightalignspace\tag 11.9
\endalign 
$$
The same inequalities also holds with ${\Cal R}$ replaced by ${\Cal U}$ everywhere 
and ${\Cal S}$ replaced by ${\Cal T}$:
$$\align 
|W|_{r}^{{\Cal U}}&\leq 
 K_1 \big(\,|\curl \underline{W}|_{r-1}^{\Cal U}+|\div W|_{r-1}^{\Cal U}
+ |W|_{r}^{\Cal T}
+\sum_{s=1}^{r}|g|_{s}^{\Cal U} |W|_{r-s}^{\Cal U}\big),
\tag 11.10\\
|W|_{r}^{{\Cal U}}&\leq 
 K_1 \sum_{s=0}^{r}\,\,
[g]^{\,\Cal U}_{s}\big(\,|\curl \underline{W}|_{r-1-s}^{\Cal U}
+|\div W|_{r-1-s}^{\Cal U}
+|W|_{r-s}^{\Cal T}\big).\rightalignspace
\tag 11.11
\endalign 
$$
\endproclaim
\demo{Proof of Lemma {11.}3 } We will first prove (11.8)
We claim that
$$
\sum_{|I|=r,U\in{\Cal R}}\!\!
|{\Cal L}_U^I W|\leq K_1\!\!\!\!\sum_{|J|=r-1, U\in{\Cal R} }
\!\!\!\big(\,|\curl\underline{{\Cal L}_U^J W}\,|
+|\div{\Cal L}_U^J W|+ [g]_1 |{\Cal L}_U^J W|\big)
+K_1\!\!\!\!\sum_{|I|=r,S\in {\Cal S}}\!\!|{\Cal L}_S^I W|\tag 11.12
$$
First we note that there is noting to prove if $d(y)\geq d_0$ since then
${\Cal S}$ span the full tangent space. Therefore, it suffices to prove (11.12)
when $d(y)\leq d_0$ and with ${\Cal S}$ replaced by ${\Cal S}_0$ and
${\Cal R}$ replaced by ${\Cal R}_0$.
Then (11.12) follows from (11.4) if $r=1$ and assuming that its true for
$r$ replaced by $r\!-\!1$ we will prove that it holds for $r$.
If we apply (11.4) to ${\Cal L}_{U}^J W$, where $|J|=r\!-\!1$, we get
$$
|\hat{\Cal L}_U {\Cal L}_U^J W|
\leq K_1\big(\,|\curl\underline{\hat{\Cal L}_U^J W}\,|
+|\div{\Cal L}_U^J W|+\sum_{S\in{\Cal S}}|{\Cal L}_S {\Cal L}_U^JW|
+[g]_1|{\Cal L}_U^J W|\big).\tag 11.13
$$
If ${\Cal L}_U^J$ consist of all tangential derivatives then
it follows that $|{\Cal L}_U \hat{\Cal L}_U^J W|$ is bounded by the right hand
side of (11.12). If ${\Cal L}_U^J$ does not consist of only
tangential derivatives then, since $[{\Cal L}_R, {\Cal L}_S]={\Cal L}_{[R,S]}=0$,
if $S\in {\Cal S}_0$, we can write 
${\Cal L}_S\hat{\Cal L}_U^J W=\hat{\Cal L}_U^K {\Cal L}_{S^\prime}W$,
for some $S^\prime\in {\Cal S}_0$.
If we now apply (11.12) with $r$ replaced by $r-1$ to ${\Cal L}_{S^\prime} W$,
(11.12) follows also for $r$.

In (11.8) we have ${\Cal L}_U^I\curl \underline{W}=\curl{\Cal L}_U^I \underline{W}$
which however is different from $\curl\,\underline{{\Cal L}_U^I W}$. 
We have:
$$
{\Cal L}_U^J \underline{W}_a={\Cal L}_U^J (g_{ab}W^b)=-g_{ab}{\Cal L}_U^J W^b
+\tilde{c}_{J_1 J_2}^{\,J}g^{J_1}_{ab}{\Cal L}_U^{J_2} W^b,
\qquad\text{where}\quad g^{J}_{ab}={\Cal L}_U^J g_{ab}\tag 11.14
$$
where the sum is over all $J_1+J_2=J$ and $\tilde{c}_{J_1 J_2}^{\,J}=1$ for $|J_2|<|J|$
$\tilde{c}_{J_1 J_2}^{\,J}=0$ if $J_2=J$.
It follows that
$$
|\curl \underline{{\Cal L}_U^J W}-\curl{\Cal L}_U^J \underline{W}\,|\leq
2\tilde{c}_{J_1 J_2}^{\, J} \big( |\pa g^{J_1}| |{\Cal L}_U^{J_2} W|
+|g^{J_1}| |\pa{\Cal L}_U^{J_2} W|\big),\qquad
|J_2|<|J|,\tag 11.15
$$
where the partial derivative can be estimated by Lie derivatives.
(11.9) follows by induction from (11.8).
Finally, (11.10) follows from (11.12) and (11.15). In fact, applying (11.12) to
$W$ replaced by ${\Cal L}_{D_t}^k W$ we see that (11.12) holds also for 
${\Cal R}$ replaced by ${\Cal U}$ and ${\Cal S}$ replaced
by ${\Cal T}$ and (11.15) also holds for $U\in {\Cal U}$. 
\enddemo

\head 12. The estimates for the curl and the normal derivatives.\endhead 
Note that in section 10 we only had bounds for the derivatives that are tangential at 
the boundary, as well as all derivatives in the interior since $S$ span the
full tangent space in the interior. We will now use estimates for the curl together 
with the estimates for the tangential derivatives to get 
estimates also for normal derivatives close to the boundary. 
Let 
$$
\dot{w}_a=g_{ab}\dot{W}^b,\quad\text{and}\quad 
\curl\, {w}_{ab}=\pa_a w_b-\pa_b w_a.\tag 12.1 
$$
Then we have 
$$
{D_t}\big(g_{ab}\dot{W}^b\big)-\pa_a\big((\pa_c p) W^c\big)
=-\pa_a q +\omega_{ab}\dot{W}^b+\underline{F}_a \tag 12.2
$$
Note that (12.2) can also be formulated as 
$$
{D_t} \dot{w}+\underline{A} W-\underline{C}\dot{W}=\underline{F}\tag 12.3
$$
where the underline as before means that we lowered the indices so
the result is a one form. Note here that $\dot{w}$ is not 
equal $D_t w$ so the notation is slightly confusing. 
But what we mean is that we think of $W$ as a vector field and take
the time derivative as a vector field which results in $\dot{W}$
and then $\dot{w}$ is the corresponding one form obtained by lowering the 
indices. 
We obtain 
$$
{D_t}\curl\dot{w}_{ab}=
(\pa_c\, {\omega}_{ab})\dot{W}^c
-\omega_{cb} \pa_a \dot{W}^c+\omega_{ca} \pa_b \dot{W}^c
+\curl \underline{F}_{ab}\tag 12.4
$$
Since ${D_t} w_a=\dot{g}_{ab}W^b+g_{ab}\dot{W}^b$ and
$\pa_a \dot{g}_{bc}-\pa_b\dot{g}_{ac}=\pa_c\omega_{ab}$ 
we also obtain 
$$
{D_t}\curl{w}_{ab}=\curl{\dot{w}}_{ab}+
(\pa_c\, {\omega}_{ab}){W}^c
+\dot{g}_{bc}\pa_a {W}^c
-\dot{g}_{ac}\pa_b {W}^c.\tag 12.5
$$
Since $\div W=\div \dot{W}=0$ it follows from Lemma {11.}2 
and (12.4)-(12.5) that 
$$\align 
|{D_t}\curl\dot{w}| &\leq K_1 |\omega|
\big( |\curl \dot{w}|+\tsize{\sum_{S\in{\Cal S}}} |{\Cal L}_S \dot{W}|
+[g]_1|\dot{W}|\big)+|\pa \omega| |\dot{W}|+|\curl \underline{F}|\tag 12.6\\
|{D_t}\curl{w}| &\leq |\curl \dot{w}|+ K_1 |\dot{g}|
\big( |\curl{w}|+\tsize{\sum_{S\in{\Cal S}}} |{\Cal L}_S {W}|
+[g]_1 |W|\big)+|\pa\omega| |W|\tag 12.7
\endalign 
$$
Since we already have control of the tangential derivatives $S$ by
section 11 
this obviously gives us control of $\curl \dot{w}$ and $\curl w$ as well
and once we have control of these we in fact control all components by 
Lemma {11.}3 again.
The norms will be measured in $L^2$ since we have control of the $L^2$ norms of
the tangential components. We will now derive higher order versions of 
the inequalities (12.6)-(12.7) using the higher order version of Lemma {11.}2,
i.e. (11.9) in  Lemma {11.}3. 

We must now get equations for the curl of higher derivatives as well. 
Applying ${\Cal L}_{U}^J$ to (12.4)-(12.5) gives, since 
the Lie derivative commutes with the curl, 
$$
{D_t}\curl{\Cal L}_{U}^J\dot{w}_{ab}=
c_{J_1 J_2}\Big((\pa_c\, {\omega}_{ab}^{J_1}) {\Cal L}_{U}^{J_2} \dot{W}^c
-\omega_{cb}^{J_1} \pa_a  {\Cal L}_{U}^{J_2} \dot{W}^c+\omega_{ca}^{J_1} 
\pa_b {\Cal L}_{U}^{J_2}\dot{W}^c\Big)
+(\curl {\Cal L}_U^J \underline{F})_{ab},\tag 12.8
$$
where $\omega^J= {\Cal L}_{U}^{J}\, \omega $ and 
$$
{D_t}\curl{\Cal L}_{U}^J{w}_{ab}=\curl{\Cal L}_{U}^J{\dot{w}}_{ab}+
c_{J_1 J_2}(\pa_c \,{\omega}_{ab}^{J_1}){\Cal L}_{U}^{J_2}{W}^c
+c_{J_1 J_2}\Big( \dot{{g}}^{J_1}_{bc}\pa_a{\Cal L}_{U}^{J_2} {W}^c
-\dot{{g}}^{J_1}_{ac}\pa_b{\Cal L}_{U}^{J_2} {W}^c\Big) \tag 12.9
$$
where $\dot{{g}}^{J}_{ab}={\Cal L}_U^{J} {D_t} g_{ab}$. Let us make a definition: 
\demo{Definition {12.}1} Let $\beta$ be a two form. 
With notation as in Definition {11.}2 we set 
$$
([g]|\beta|)_{u}^{\Cal V} =\sum_{s+r\leq u}
 [g]_{s}^{\Cal V}|\beta|_{r}^{\Cal V}.\tag 12.10
$$
\enddemo
Using Lemma {11.}3 and Lemma {11.}2 it follows that:
\proclaim{Lemma {12.}1} With notation as in Definition {11.}1
and  Definition {12.}1 and the convention that 
$|\curl \underline{W}|_{-1}^{\Cal V}=|\div W|_{-1}^{\Cal V}=0$ we have 
$$\align 
\leftalignspace |D_t \curl \dot{w}|_{r-1}^{\Cal R}&
\leq  K_1 \sum_{s=0}^{r}\,\,
([g]|\omega|)^{\Cal R}_{r-s}\big(|\curl \dot{w}|_{s-1}^{\Cal R}+|\div
\dot{W}|_{s-1}^{\Cal R}+
|\dot{W}|_{s}^{\Cal S}\big)+|\curl \underline{F}|_{r-1}^{\Cal R}
\rightalignspace \tag 12.11\\
\leftalignspace |D_t \curl w|_{r-1}^{\Cal R}&\leq K_1 \sum_{s=0}^{r}
\,\,
([g]|\dot{g}|)^{\Cal R}_{r-s}\big(|\curl w|_{s-1}^{\Cal R}+|\div W|_{s-1}^{\Cal R}+
|W|_{s}^{\Cal S}\big)+ |\curl \dot{w}|_{r-1}^{\Cal R}
\rightalignspace \tag 12.12
\endalign 
$$
The same inequalities hold with ${\Cal R}$ replaced by ${\Cal U}$ and 
${\Cal S}$ replaced by ${\Cal T}$. 
\endproclaim 
\demo{Proof of Lemma {12.}1} Let us first prove (12.11).
The first terms in the right hand side of (12.8)
are by Lemma {11.}3  bounded by a constant times 
$$
\sum_{u=0}^r |\omega|_{r-u}^{{\Cal R}} |\dot{W}|_{u}^{{\Cal R}}\leq 
K_1 \sum_{u=0}^r \sum_{s=0}^{u}
|\omega|_{r-u}^{{\Cal R}} 
[g]^{\Cal R}_{u-s}\big(|\curl \dot{w}|_{s-1}^{\Cal R}+|\div \dot{W}|_{s-1}^{\Cal R}+
|\dot{W}|_{s}^{\Cal S}\big)\tag 12.13
$$
The proof of (12.12) uses the same argument and that 
that $\pa_c\omega_{ab}=\pa_a \dot{g}_{bc}-\pa_b\dot{g}_{ac}$. 
\qed\enddemo 
Let us now introduce some new norms and some new notation:

\demo{Definition {12.}2} For ${\Cal V}$ any of our families of vector fields let 
$$
\|W\|^{{\Cal V}}_{\,r}=\|W(t)\|_{{\Cal V}^{\,r}(\Omega)}=\sum_{|I|\leq r, I\in {\Cal V}}
\Big({\int_\Omega |{\Cal L}_U^I W(t,y)|^2\, \kappa dy}\Big)^{1/2},\tag 12.14
$$
and 
$$
C_r^{\Cal{V}}=\sum_{|J|=\leq  r-1,\, J\in {\Cal V}}\Big(\int_\Omega |\curl{\Cal L}_U^J  \dot{w}|^2
+|\curl{\Cal L}_U^J {w}|^2\, \kappa dy\Big)^{1/2},
\qquad C_0^{\Cal V}=0.
\tag 12.15
$$
\enddemo 
Note that $\|W(t)\|_{{\Cal R}^{\,r}(\Omega)}$
is equivalent to the usual Sobolev norm in the Lagrangian coordinates.
\demo{Definition {12.}3} For ${\Cal V}$ any of our families of vector fields 
and for $\beta$ a function, a 1-form, a 2-form or a vector field let $|\beta|_s^{\Cal V}$
be as in Definition {11.}1 and set 
$$
\|\beta\|_{s,\infty}^{\Cal V}=\| \, |\beta|_s^{\Cal V} \|_{L^\infty(\Omega)},
\qquad 
[[g]]_{s,\infty}^{\Cal V}=\!\!\!\!\!\!\!\!\sum_{s_1+...+s_{k}\leq s,
\, s_i\geq 1}\!\!\!\!\!\!\!\!
 \|g\|_{s_1,\infty}^{\Cal V}\!\!\cdot\cdot\cdot \|g\|_{s_{k},\infty}^{{\Cal
V}}, 
\qquad [[g]]_{0,\infty}^{\Cal V}=1\tag 12.16
$$
where the sum is over all combinations with $s_i\geq 1$.  
Furthermore, let 
$$
m_{r}^{\Cal V}=[[g]]_{r,\infty}^{\Cal V},
\qquad\quad \dot{m}_{r}^{\Cal V}=\sum_{s+u\leq r}[[g]]_{s,\infty}^{\Cal V}
\big(||\dot{g}||_{u,\infty}^{\Cal V}+||\omega||_{u,\infty}^{\Cal V}\big). 
\tag 12.17 
$$
\enddemo 
Let 
$F_r^{\Cal U}=\| \curl \underline{F}\|_{{\Cal U}^{r-1}(\Omega)}$. 
It now follows from Lemma {12.}1 that 
$$
\Big|\frac{d C^{\,\Cal{U}}_r}{dt}\Big|\leq 
K_1 \sum_{s=0}^{r} \dot{m}^{\Cal U}_{r-s} \big(C^{\,\Cal{U}}_s
+ E^{\Cal T}_s\big)
+F_r^{\Cal U}\tag 12.18
$$
where $E^{\Cal T}_s$ is the energy of the tangential derivatives defined in section 10. 
Hence 
$$
C^{\Cal{U}}_r\leq K_1 e^{\int_0^t K_1 \, \dot{m}^{\Cal U}_0 \, d\tau} 
 \int_0^t \Big(\sum_{s=1}^{r-1} \dot{m}^{\Cal U}_{r-s} C^{\Cal{U}}_s
+\sum_{s=0}^r  \dot{m}^{\Cal U}_{r-s} E^{\Cal T}_s
+F_r^{\Cal U}\Big)\, d\tau \tag 12.19 
$$
Since we already proved a bound for $E^{\Cal T}_s$ in Lemma {10.}1 it 
inductively follows that $C_r^{\Cal U}$ is bounded. 
Note that, if $r=1$ the interpretation of (12.19) is that the first sum is not there.
By Lemma {11.}3:
$$
\| W(t)\|_{{\Cal U}^r(\Omega)}+\| \dot{W}(t)\|_{{\Cal U}^r(\Omega)}\leq 
K_1 \sum_{s=0}^{r} m^{\Cal U}_{r-s}\big(C^{\Cal{U}}_s+E_{s}^{\Cal T}\big)
\tag 12.20
$$
Hence we have: 
\proclaim{Lemma {12.}2}  Suppose that 
$x, p\in C^{r+2}([0,T]\times\Omega)$, $p\big|_{\pa\Omega}=0$,
$\na_N p\big|_{\pa\Omega}\leq -c_0<0$ and $\div V=0$, where $V=D_t x$.
Then there is a constant $C=C(x,p)$ depending only on the norm 
of $(x,p)$, a lower bound for $c_0$ and an upper bound for $T$, 
such that if $E^{\Cal T}_s(0)=C^{\Cal U}_s(0)=0$, for $s\leq r$, then 
$$
 \| W\|_r^{\Cal U}+
\| \dot{W}\|_r^{\Cal U}+E^{\Cal T}_r
\leq C\int_0^t   \|F\|_{r}^{\Cal U}\, d\tau 
,\qquad\text{for}\quad 0\leq t\leq T. 
\tag 12.21 
$$
\endproclaim

\head 13. The smoothed out normal operator. \endhead 
In order to prove existence we first have to replace the normal operator $A$
by a sequence $A^\varepsilon$ of bounded symmetric and positive operators
that convergence to $A$, as $\varepsilon\to 0$.
The boundedness is needed for the existence and
the symmetry and positivity is needed to get a positive term in the energy.
Furthermore the commutators with Lie derivatives with respect to tangential 
vector fields as well as the curl
have to be well behaved. 
Let $\rho=\rho(d)$ be a smooth function of $d=d(y)=\dist{(y,\pa\Omega)}$, 
such that 
$$
\rho^\prime\geq 0,\qquad \rho(d)=d\quad\text{for $d\leq 1/4$}
\quad \text{and}\quad \rho(d)=1/2\quad
\text{for $d\geq 3/4$}.\tag 13.1
$$
 Let $\chi(\rho)$ be a smooth function such that 
$$
\chi^\prime(\rho)\geq 0,\qquad \chi(\rho)=0, \quad
\text{when}\quad \rho\leq
1/4,\quad\text{and}\quad 
\chi(\rho)=1,\quad\text{when}\quad \rho\geq 3/4\tag 13.2
$$
For a function $f$ vanishing on the boundary we define
$$
A_f^\varepsilon W^a=
P\big(-g^{ab}\chi_\varepsilon(\rho)\pa_b \big( f \rho^{-1}
 (\pa_c \rho) W^c\big)\big)\tag 13.3
$$
where $\chi_\varepsilon(\rho)=\chi(\rho/\varepsilon)$. 
Then if we integrate by parts we get 
$$
\langle U,A_f^\varepsilon W\rangle =\int_\Omega f \rho^{-1} \chi^\prime_\varepsilon(\rho)
(U^a\pa_a \rho)(W^b\pa_b \rho) dy \tag 13.4
$$
from which it follows that $A_f^\varepsilon$ is symmetric and 
$$
A^\varepsilon=A_p^\varepsilon\geq 0,\qquad\text{i.e.}\qquad
\langle W,A^\varepsilon W\rangle \geq 0,\qquad\text{if}\quad p\geq 0\tag 13.5
$$
It also follows that another expression for $A_f^\varepsilon$ is
$$
A_f^\varepsilon W^a=
P\big(g^{ab}\chi_\varepsilon^\prime(\rho)(\pa_b\rho) 
f \rho^{-1} (\pa_c \rho) W^c\big)\tag 13.6
$$

$A^\varepsilon $ is now for each $\varepsilon>0$ a bounded operator
$$
\| A^\varepsilon W\|\leq C\|\na_N p\|_{L^\infty} \varepsilon^{-1}\|W\|\tag 13.7
$$
since $\chi_\varepsilon^\prime\leq C/\varepsilon$ and
$ p\, \rho^{-1} |\pa \rho|\leq C \|\na_N p\|_{L^\infty(\pa \Omega)}$. 
In general, since the projection is continuous on $H^r(\Omega)$, 
see (3.6) and (3.8), 
if the metric and pressure are sufficiently regular we get 
$$
\sum_{j=0}^k\|D_t^j  A^\varepsilon W\|_{H^r(\Omega)}\leq 
C_{\varepsilon,r,k} \sum_{j=0}^k\|D_t^j W\|_{H^r(\Omega)}. \tag 13.8
$$
Moreover 
$$
A^\varepsilon U\to AU,\quad\text{in}\quad L^2(\Omega),\quad
\text{if}\quad U\in H^1(\Omega)\tag 13.9
$$
In fact, the projection is continuous in the norm and 
$\chi_\varepsilon F\to F$ in $L^2$ if $F\in L^2$. 
It follows that 
$$
P\big(g^{ab}\chi_\varepsilon(\rho)\pa_b \big( 
 p\, \rho^{-1} (\pa_c \rho) U^c\big)\big)\to 
P\big( g^{ab}\pa_b \big( p\, \rho^{-1}
 (\pa_c \rho) U^c\big)\big)
=P\big(g^{ab}\pa_b \big( (\pa_c p) U^c\big)\big)\tag 13.10
$$
since $p \, \rho^{-1} \pa_c \rho=\pa_c p$ on the boundary. 

We will now calculate the commutators with the Lie derivative ${\Cal L}_T$ with 
respect to tangential vector fields $T$.
As before the inequality 
$$
|\langle U, A_{fp}^\varepsilon W\rangle |\leq 
\langle U, A_{|f|p}^\varepsilon U\rangle^{1/2}
\langle W, A_{|f|p}^\varepsilon W\rangle^{1/2}
\leq \|f\|_{L^\infty(\Omega\setminus\Omega^\varepsilon)} 
\langle U, A^\varepsilon U\rangle^{1/2}\langle W, A^\varepsilon W\rangle^{1/2}\tag
13.11
$$
hold, where $\Omega^\varepsilon=\{y\in\Omega;\, d(y)>\varepsilon\}$.
In fact, it suffices to take the supremum over the set where $d(y)\leq 
\varepsilon$ since $\chi^\prime_\varepsilon=0$, when $d(y)\geq \varepsilon$. 
  The only difference with (3.13) is that now the supremum over 
a small neighborhood of the boundary instead of on the boundary.
The positivity properties (13.5) and (13.11) for $A^\varepsilon$
will play the role that
(3.12) and (3.13) did for $A$.
In particular, since $p$ vanishes on the boundary,
$p>0$ in the interior and $\na_N p\leq -c_0<0$ on the boundary it follows that 
$\dot{p}=D_t\, p$ vanishes on the boundary and $\dot{p}/p$ is a smooth function. 
Therefore 
$$
\dot{A}^\varepsilon=A_{\dot{p}}^\varepsilon\qquad\text{satisfies}\qquad 
|\langle W, \dot{A}^\varepsilon W\rangle |
\leq \|\dot{p}/ p\|_{L^\infty(\Omega\setminus\Omega^\varepsilon)}
\langle W, A^\varepsilon W\rangle \tag 13.12
$$
Here 
$\dot{A}^\varepsilon$ is the time derivative of the operator $A^\varepsilon$, 
considered as an operator with values in the one forms.
It will show up in the energy estimate for the $\varepsilon$ smoothed out equation 
in the next section. 

The commutators between $A_f^\varepsilon$ and Lie derivatives with 
respect to tangential vector fields are basically the same as for $A$.
Note that 
$$
T d=0,\quad 
\text{if} \quad T\in {\Cal T}_0={\Cal S}_0\cup\{D_t\}\tag 13.13
$$
where ${\Cal S}_0$ are the rotations. Hence if $T$ is any of these vector fields 
we have 
$$
P\big(g^{ca}{\Cal L}_{T}(g_{ab} A_f^\varepsilon W^b)\big)
=A_f^\varepsilon {\Cal L}_T W^c + A^\varepsilon _{Tf} W^c,\tag 13.14
$$
However, in order to get additional 
regularity in the interior we include the vector fields ${\Cal S}_1$ 
that span the tangent space in the interior.
The vector fields in ${\Cal S}_1$ satisfy 
$$
S\rho ={\Cal L}_S \rho=0,\quad \text{when}\quad d\leq d_0/2\tag 13.15
$$
Since $\chi^\prime_\varepsilon(\rho)=0$ when $d\geq \varepsilon$
 the commutator relation (13.14) above is true for these as well if we assume that
$\varepsilon\leq d_0/2$.

It remains to estimate the curl of $A^\varepsilon$.
Whereas, the curl of $A$ vanishes this is not the case for the 
curl of $A^\varepsilon$. It will however vanish away from the boundary.
With $\underline{A}^\varepsilon W_a=g_{ab} A^\varepsilon W^b$ we have 
$$
\underline{A}^\varepsilon W_a=-\chi_\varepsilon(\rho)\pa_a \big( 
  p\, \rho^{-1} (\pa_c \rho) W^c\big)-\pa_a q_1\tag 13.16
$$
for some function $q_1$ vanishing on the boundary and determined
so the divergence vanishes. Since the curl of the gradient vanishes and
$\chi^\prime_\varepsilon(\rho)=0$ when $d\geq \varepsilon$ we have
$$
\curl \, \underline{A}^\varepsilon W_{ab} =0,\qquad\text{when}\quad d(y)\geq \varepsilon
\tag 13.17
$$

\head 14. The smoothed out equation and existence of weak solutions.\endhead 
The $\varepsilon $ smoothed out linear equation:
$$
\ddot{W}_\varepsilon^a+A^\varepsilon W_\varepsilon^a
+\dot{G}\dot{W}_\varepsilon-C\dot{W}_\varepsilon^a=F^a,
\qquad \dot{W}_\varepsilon\Big|_{t=0}={W}_\varepsilon\Big|_{t=0}=0\tag 14.1
$$
is just an ordinary differential equation for $(W_\varepsilon,\dot{W}_\varepsilon)$
on the space of divergence free 
vector fields in $L^2(\Omega)$ since all operator 
are bounded so existence follows in $L^2(\Omega)$. 
In fact its an ordinary differential equation in the Sobolev spaces 
$H^r(\Omega)$ by (13.8).
To get additional regularity in time as well we apply more 
time derivatives using (13.8) and (3.8) and that
the initial conditions for these vanishes as well
since we constructed $F$ in (14.1) so it vanishes to any given order.
If initial data, encoded in $F$, are smooth, we hence have a smooth
solution of the $\varepsilon$ approximate linear equation. 

Now we want to use the existence and estimates
 for the $\varepsilon$ smoothed out linear equation
and pass to the limit as $\varepsilon \to 0$ to get existence
for the linearized equation. Will show that $W_\varepsilon\to W$
weakly in $L^2$, where $W\in H^r(\Omega)$ for some large $r$.
From the weak convergence it will follow that $W$ is a weak solution
and then from the additional regularity of $W$ it will follow that in 
fact its a classical solution and hence that the {\it a priori } bounds in 
the earlier section hold.

Of course the norm
of $A^\varepsilon$ tends to infinity as $\varepsilon\to 0$
but since it is a positive operator it can be included in the energy.
The energy will be the same as before
with $A$ replaced by $A^\varepsilon$, so (4.4) becomes
$$
E^\varepsilon=\langle \dot{W}_\varepsilon,\dot{W}_\varepsilon\rangle
 +\langle W_\varepsilon, (A^\varepsilon+I) W_\varepsilon\rangle \tag 14.2
$$
The time derivative of the first term is the same as (4.5) with $W$ replaced 
by $W_\varepsilon$. 
Since $D_t\, d=0$ it follows from taking the time derivative of (13.4),
with $f=p$, that 
$$
\frac{d }{dt} \langle W_\varepsilon, A^\varepsilon W_\varepsilon\rangle
=2 \langle \dot{W}_\varepsilon, A^\varepsilon W_\varepsilon\rangle
+ \langle W_\varepsilon, A^\varepsilon_{\dot{p}} W_\varepsilon\rangle,\tag 14.3
$$
where the last term is bounded by (13.12). Hence by (4.7)-(4.9): 
$$
|\dot{E}^\varepsilon|\leq \Big(1+\|\dot{g}\|_{L^\infty(\Omega)}
+\|( D_t p)/ p\|_{L^\infty(\Omega)}\Big) E^\varepsilon+2\sqrt{E^\varepsilon}
\|F\|\tag 14.4
$$
from which we get a uniform bound for $0\leq t\leq T$ independent of 
$ \varepsilon$: $E^\varepsilon(t)\leq C $.

Since $\|W_{\varepsilon}\| \leq C$ we can now choose a 
subsequence $W_{\varepsilon_n}\to W$ weakly in the inner product. 
We will show below that the limit $W$ is a weak solution if the equation. 
Multiplying the $\varepsilon$ smoothed out equation 
by a smooth divergence free vector field $U$ that vanishes 
for $t\geq T$ and integrating by parts we get
$$
\int_0^T\int_{\Omega}  g_{ab}\Big(\ddot{U}^a+A^\varepsilon U^a
+\dot{G} \dot{U}^a-C\dot{U}^a-\dot{C} U^a \Big) W_\varepsilon^b\, dy dt
=\int_0^T\int_{\Omega} g_{ab} U^a F^b\, dy dt\tag 14.5
$$
where $\dot{C} W^c=P\big(g^{ac}\dot{\omega}_{cb} W^b)$, 
since $A^\varepsilon$ and $D_t^2 +B{D_t}$
are symmetric and the adjoint of $C{D_t}$ is 
$C{D_t}-\dot{C}$. 
We proved in the previous section that 
 $A^\varepsilon U$ converges to $A U$ strongly in the 
norm if $U$ is in $H^1$. Since $W_{\varepsilon_n}\to W$ weakly this proves that 
we have a weak solution $W$ of the equation: 
$$
\int_0^T\int_\Omega g_{ab}\Big(\ddot{U}^a+A U^a
+\dot{G} \dot{U}^a-C\dot{U}^a-\dot{C} U^a \Big) W^b\, dy dt=
\int_0^T\int_\Omega g_{ab} U^a F^b\, dy dt\tag 14.6
$$
for any divergence free smooth vector field $U$ that vanishes for $t\geq T$.
Furthermore since $W_\varepsilon$ is divergence free, we have 
$$
\int_0^T\int_\Omega (\pa_a q) W_\varepsilon^a\, dy dt=0\tag 14.7
$$
for any smooth $q$ that vanishes on the boundary and hence 
$$
\int_0^T\int_\Omega (\pa_a q) W^a\, dy dt=0\tag 14.8
$$
so $W$ is weakly divergence free. 

\head 15. Existence of smooth solutions for the linearized equation.\endhead 
Now that we have existence of a weak solution we will prove that we have additional
regularity and in fact that, $W,\dot{W}\in H^{r}(\Omega)$ for any  $r\geq 0$.
It then follows that we can integrate by parts again in the above 
integrals and conclude that
$$
\int_0^T \int_\Omega q \, \pa_a W^a\, dy dt=0\tag 15.1
$$
for any smooth function $q$ that vanishes on the boundary. 
Hence $W$ is divergence free. Furthermore 
$$
\int_0^T\int_\Omega g_{ab}
U^a \Big(\ddot{W}^b+A W^b+\dot{G}\dot{W}^b-C\dot{W}^b\Big)
 \, dy dt=\int_0^T\int_\Omega g_{ab} U^a F^b\, dy dt\tag 15.2
$$
for any smooth divergence free vector field $U$ that vanishes for $t\geq T$.
But in fact since $W$ is divergence free it follows that 
$\ddot{W}^b+A W^b+\dot{G}\dot{W}^b-C\dot{W}^b$ is divergence 
and since by construction $F$ is divergence free as well it follows that (15.2)
holds for any smooth vector field $U$ that vanishes for $t\geq T$.
We then conclude that 
$$
\ddot{W}^b+A W^b+\dot{G}\dot{W}^b-C\dot{W}^b=F^b,\qquad\quad 
 \div W=0\tag 15.3
$$

It therefore only remains to show that $W\in H^r(\Omega)$.
We must show that we have uniform bounds for 
the $\varepsilon$ smooth out equation similar to the {\it a priori} bounds 
for the linearized equation.

The uniform tangential bounds for the $\varepsilon$ smoothed out equation 
follows the proof of the {\it a priori} tangential bounds in section 10. 
The proof is just a change of notation. Let 
$$
E_I^{\varepsilon}= \langle \dot{W}_{\varepsilon I},\dot{W}_{\varepsilon I}\rangle +
\langle W_{\varepsilon I}, (A+I)W_{\varepsilon I}\rangle,
\qquad W_{\varepsilon I}={\Cal L}_T^I W_\varepsilon.\tag 15.4 
$$
If $\varepsilon<d_0$ then the commutator relation for $A^\varepsilon$,
(13.14),  is exactly the same as for
$A$, (8.5). Furthermore the positivity property 
for $A_f^\varepsilon$ only differs from the one for $A_f$ 
by that the supremum over the boundary in (3.13) is replaced 
by the supremum over 
a neighborhood of the boundary where $d(y)<\varepsilon$ in (13.11). 
Hence all the calculations and inequalities in sections 10 and 12 hold
with $A$ replaced by $A^\varepsilon$, if we replace the supremum of 
$\na_N q/\na_N p$ over the 
boundary in (10.9) by the supremum of $q/p$ over the domain $\Omega\setminus\Omega^\varepsilon$,
where $\Omega^\varepsilon$ is given by (15.6). Therefore we will arrive 
at the energy bound (10.21) for $E^{{\Cal T}}_r$ replaced by
$$
E^{{\Cal T},\varepsilon}_r=
\sum_{|I|\leq r,\, I\in{\Cal T}} \sqrt{E_I^{\varepsilon}},\tag 15.5 
$$
i.e. Lemma {10.}1 hold for $E^{{\Cal T}}_r$ replaced by $E^{{\Cal T},\varepsilon}_r$
with a constant independent of $\varepsilon$. 
Note that, this is where we 
need to have vanishing initial conditions and an inhomogeneous 
term that vanishes to high order when $t=0$ so that also the higher 
order time derivatives of the solution of (14.1) vanished when $t=0$. 
If the initial conditions for higher order time derivatives were
to be obtained from the $\varepsilon$ smoothed out equation, then
they would depend on $\varepsilon$ and so we would not have been able to 
get a uniform bound for the energy, $E^{{\Cal T},\,\varepsilon}_r$. 

The bound for curl is very simple since by (13.17) the curl of $A_\varepsilon$ 
vanishes in 
$$
\Omega^\varepsilon=\{y;\, \operatorname{dist}(y,\pa\Omega)>\varepsilon\},\tag 15.6
$$
it follows that all the formulas in section 12 hold when $d(y)\geq \varepsilon$. 
This follows from
replacing $\underline{A}$ in (12.3) by $\underline{A}^\varepsilon$
and using that the curl of this vanishes for $d(y)\geq \varepsilon$.
Since all the estimates used from section 11 are point wise estimates 
we conclude that (12.11)-(12.12) hold for $W$ replaced by $W^\varepsilon$ 
when $d(y)\geq \varepsilon$. 
Let 
$$
C^{\Cal{U},\,\varepsilon}_r=
\sum_{|J|\leq r-1,\, J\in{\Cal U}}\Big(\int_{\Omega^\varepsilon} |\curl {\Cal L}_U^J
{w_\varepsilon}|^2+
 |\curl {\Cal L}_U^J \dot{w}_\varepsilon|^2\, dy\Big)^{1/2}\tag 15.7
$$
With  $C_r^{\Cal{U}}$ replaced by $C_r^{\Cal{U},\,\varepsilon}$ 
and $E^{\Cal T}_s$ replaced by $E^{\Cal T,\, \varepsilon }_s$
we get exactly the same inequalities as before (12.18)-(12.19),
since these were derived from the point wise bounds in section 11. 
Furthermore, the inequality (12.20) hold 
as well if we replace the norms by 
$$
\|W(t)\|_{{\Cal U}^{\,r}(\Omega^{\varepsilon})}=\sum_{|I|\leq r,\, I\in {\Cal U}}
\Big({\int_{\Omega^{\varepsilon}}
 |{\Cal L}_U^I W(t,y)|^2\, \kappa dy}\Big)^{1/2},\tag 15.8
$$
Therefore we conclude that the inequality in Lemma {12.}2 hold
with a constant $C$ independent of $\varepsilon$ if we replace the norms by (15.8): 
\proclaim{Lemma {15.}1}  Suppose that 
$x, p\in C^{r+2}([0,T]\times\Omega)$, $p\big|_{\pa\Omega}=0$,
$\na_N p\big|_{\pa\Omega}\leq -c_0<0$ and $\div V=0$, where $V=D_t x$.
Suppose that $W_\varepsilon$ is a solution of (14.1) where $F$ is divergence free and 
vanishing to order $r$ as $t\to 0$. Let $E^{\Cal T,\,\varepsilon}_s$ be defined by (15.5).
Then there is a constant $C$ depending only on the norm 
of $(x,p)$, a lower bound for $c_0$ and an upper bound for $T$,
but independent of $\varepsilon$, 
such that if $E^{\Cal T,\varepsilon}_s(0)=C^{\Cal U,\varepsilon}_s(0)=0$, for $s\leq r$, then 
$$
\| W_\varepsilon(t)\|_{{\Cal U}^r(\Omega^{\varepsilon})}+
\| \dot{W}_\varepsilon(t)\|_{{\Cal U}^r(\Omega^{\varepsilon})}+E^{{\Cal
T},\, \varepsilon}_r(t)
\leq C\int_0^t   \|F\|_{r}^{\Cal U}\, d\tau
 ,\qquad\text{for}\quad 0\leq t\leq T. 
\tag 15.9
$$
\endproclaim

It therefore follows that the limit $W$ satisfies the same bound with 
$\Omega^\varepsilon$ replaced by $\Omega$, and so the weak solution
in section 14 is in fact a smooth solution. 

\head 16. The energy estimate revisited and the proof of the theorem.\endhead
In section 10 we estimated the energies of the tangential derivatives
without using the estimate of the normal derivatives coming from the curl.
This was necessary to get uniform bounds 
for the $\varepsilon$ smoothed out
equation since in that case we could not estimate the curl close to the boundary.
The drawback was that instead we had to include all time derivatives as well
in the energy. However, now that we have existence we can obtain other bounds
for the linearized equation directly.
In section 9 we calculated the commutator between the linearized operator, 
considered as an operator from the divergence free vector fields to the one forms,
and Lie derivatives with respect to tangential vector fields,
and then projected the result back onto the divergence free vector fields. 
This was needed because the commutator between Lie derivatives and the 
operator $A$ considered as an operator with values in the one forms
is better behaved. 
However, the drawback is that the commutator with the second time derivative,
considered as an operator with values in the one forms,
involves second time derivatives, which is why we had to include all the time
derivatives. Now we will instead commute through directly with the operator from
the divergence free vector fields to the divergence free vector fields.
Let us then also consider the original setting with non vanishing 
initial conditions and an inhomogeneous term:
$$
 \ddot{W}^a -g^{ab}\pa_b\big( (\pa_c p) W^c -q_1\big)
=-g^{ab}\big((\dot{g}_{cb}-\omega_{cb}) \dot{W}^c-\pa_b q_2\big)+F^a\tag 16.1 
$$
where $q_1$ and $q_2$ vanishes on the boundary and are chosen so that 
each term is divergence free. The second term on the left is 
$AW^a$ and the term in the right is $-\dot{G} \dot{W}^a+C\dot{W}^a$. 
Let us now first calculate the commutators with $A$ and tangential vector fields.
$$\multline 
{\Cal L}_{S} \Big( g^{ab} \pa_b\big( (\pa_c p) W^c -q_1\big)\Big)\\
=({\Cal L}_{S} g^{ab} )\pa_b\big( (\pa_c p) W^c -q_1\big)
+ g^{ab}\pa_b\big( (\pa_c S p) W^c 
+(\pa_c p)({\Cal L}_{S} W)^c - S q_1\big)
\endmultline \tag 16.2 
$$
where ${\Cal L}_S g^{ab}\!=\!-g^{ac}g^{bd}{g}^S_{cd},\,\,$ ${g}^S_{cd}\!=\!{\Cal L}_S g_{cd}$.
Projecting each term onto divergence free vector fields: 
$$
{\Cal L}_S AW=-G_S AW+A_S W+ A{\Cal L}_S W,\tag 16.3 
$$
where $A_S=A_{Sp}$ and $G_S=M_{{g}^S}$ is the operator 
$G_S W^a=P\big(g^{ac} {g}^S_{cb} W^b\big)$. Expressed differently 
$$
[{\Cal L}_S, A ] W=(A_{S}-G_S A) W\tag 16.4
$$
Although $G_S$ is a bounded operator, all the  
positivity properties of $A$ are lost and the best we can say is that 
$G_S A$ is an operator of order $1$. The operator $A_{S}$ is also of order $1$ but 
in section 10 we used the positivity property to estimate it in terms
of $A$ which we controlled by the energy. 
It remains to calculate the commutator with 
$G_S$ and $C$, which basically  are the same. 
$$\multline 
{\Cal L}_T G_S W^i={\Cal L}_T \Big(g^{ab}\big( {g}^S_{bc} W^c-\pa_b q\big)\Big)\\
=({\Cal L}_T g^{ab})\big( {g}^S_{bc} W^k-\pa_b q\big)
+g^{ab}({\Cal L}_T{g}^S_{bc}) W^c
+g^{ab}{g}^S_{bc} {\Cal L}_T W^c-g^{ab} \pa_b T q\\
\endmultline \tag 16.5 
$$
Projecting each term onto the divergence free vector fields we arrive at 
$$
[{\Cal L}_T, G_S] W=(G_{TS}-G_T G_S) W,\tag 16.6 
$$
where $G_{TS} W^a =P\big( g^{ab}{g}^{TS}_{bc} W^c\big)$ and 
${g}^{TS}_{bc} ={\Cal L}_T{\Cal L}_S g_{bc}$. 

In general using (16.4) and (16.6) to 
commute through we get for some constants $\tilde{d}^{\,I_1...I_k}_I$
$$
{\Cal L}_S^I A W-A{\Cal L}_S^I  W=\tilde{d}^{\, I_1 I_k}_I G_{I_1}\cdot\cdot \cdot 
G_{I_{k-2}} A_{I_{k-1}} W_{I_{k}}\tag 16.7
$$
where the sum is over all combinations with $I_1+...+I_{k}=I$, with $k\geq 2$,
and $|I_k|<|I|$. 
Here $G_{J}W^a=M_{{g}^J}W^a=P(g^{ac}{g}_{cb}^J W^b)$, where 
${g}^J_{ac}={\Cal L}_S^J g_{ac}$, $A_J=A_{S^J p}$
and $W_J={\Cal L}_S^J W$. Similarly we get the commutators with $\dot{G}$ and 
$C$ 
$$
{\Cal L}_S^I \dot{G} W-\dot{G}
{\Cal L}_S^I  W=\tilde{e}^{\, I_1 I_k}_I G_{I_1}\cdot\cdot \cdot 
G_{I_{k-2}} \dot{G}_{I_{k-1}} W_{I_{k}}\tag 16.8
$$
$$
{\Cal L}_S^I C W-C
{\Cal L}_S^I  W=\tilde{e}^{\, I_1 I_k}_I G_{I_1}\cdot\cdot \cdot 
G_{I_{k-2}} C_{I_{k-1}} W_{I_{k}}\tag 16.9 
$$
The only thing that matters is that these are bounded operators,
and in fact they are lower order since $|I_k|<|I|$. 
Hence we obtain
$$
L_1 W=\ddot{W}_I+A W_I+\dot{G} \dot{W}_I-C \dot{W}_I=H_I\tag 16.10 
$$
where 
$$
\multline 
H_I=F_I
+\tilde{d}^{\, I_1 I_k}_I G_{I_1}\!\!\cdot\cdot \cdot 
G_{I_{k-2}} A_{I_{k-1}} W_{I_{k}}\\
+\tilde{e}^{\, I_1 I_k}_I G_{I_1}\!\!\cdot\cdot \cdot 
G_{I_{k-2}} \dot{G}_{I_{k-1}} \dot{W}_{I_{k}}
+\tilde{e}^{\, I_1 I_k}_I G_{I_1}\!\!\cdot\cdot \cdot 
G_{I_{k-2}} C_{I_{k-1}} \dot{W}_{I_{k}}
\endmultline \tag 16.11
$$
where $|I_k|<|I|$ in the right hand side. 
Here $F_I={\Cal L}_T^I F$. 
As before, let 
$$
E_I=\langle \dot{W}_I,\dot{W}_I\rangle +\langle W_I, (A+I)W_I\rangle \tag 16.12
$$
where we now only consider $W_I={\Cal L}_S^I W$ with $S\in {\Cal S}$.
The energy estimate is like before and we only have to be able to estimate the 
$L^2$ norm of the right hand side of (16.10). The terms on the second
row of (16.11) 
are obviously bounded by $E_J$ for some $|J|\leq |I|$. In fact they are even
lower order since we have strict inequality. Therefore it only remains to
estimate the term on the right in the first row. $|I_k|<|I|$ but $A_{I_k}$ is 
order one and it contains derivatives in any direction 
so that term has to be estimated by the $\|\pa W_{I_k}\|_{L^2(\Omega)}$,
and so it does not directly help to have an estimate for 
$\|{\Cal L}_S W_{I_k}\|_{L^2(\Omega)}$ for all tangential derivatives $S$. 
However the estimate of the tangential derivatives
 together with the estimates for curl in Lemma {12.}1  
gives the required estimate. 

Let $C_r^{\Cal R}$ be defined 
(12.15), let $E_s^{\Cal S}$ be defined by (10.13) and let $m^{\Cal R}_s$
and $\dot{m}^{\Cal R}_s$
 be as in Definition {12.}3. 
Then by Lemma {11.}3 we get the inequality corresponding to (12.20):
$$
\| W\|_{r}+\| \dot{W}\|_{r}\leq 
K_1 \sum_{s=0}^{r} m^{\Cal R}_{r-s}\big(C^{\Cal{R}}_s+E_{s}^{\Cal S}\big),
\qquad \text{where}\quad \|W\|_r= \| W(t)\|_{{\Cal R}^r(\Omega)}
\tag 16.13
$$ 
Since the projection has norm $1$, 
$\| G_{J} W\|\leq \|g^J\|_{\infty}\|W\|$. 
It follows that
$$\align 
\leftalignspace \| G_{I_1}\cdot\cdot \cdot 
G_{I_{k-2}} C_{I_{k-1}} \dot{W}_{I_{k}}\|
&\leq \| g^{I_1}\|_{\infty}\cdot\cdot\cdot
\|g^{I_{k-2}}\|_{\infty}\,\|\omega^{I_{k-1}}\|_{\infty}\,
 \|\dot{W}_{I_{k}} \|\leq  \dot{m}_{r-s}^{\Cal R} \|\dot{W}\|_{s}\rightalignspace 
\tag 16.14 \\
\leftalignspace\! \| G_{I_1}\cdot\cdot \cdot 
G_{I_{k-2}} \dot{G}_{I_{k-1}} \dot{W}_{I_{k}}\|
&\leq \| g^{I_1}\|_{\infty}\cdot\cdot\cdot
\|g^{I_{k-2}}\|_{\infty}\,\|\dot{g}^{I_{k-1}}\|_{\infty}\,
 \|\dot{W}_{I_{k}} \|
\leq  \dot{m}_{r-s}^{\Cal R} \|\dot{W}\|_{s}\!\rightalignspace\tag 16.15 
\endalign 
$$
where $s=|I_k|<r$ and $r=|I|$. 
Let
$$
p_{r}^{\Cal R}=\sum_{s=0}^r\,\, [[g]]_{r-s,\infty}^{\Cal R} 
\sum_{|J|\leq s+1,\, J\in {\Cal S}}\|\pa S^J p\|_{L^\infty(\pa\Omega)}
\tag 16.16
$$
Since $A_{J}=A_{S^J p}$ it follows form (3.15) that 
$$
\| G_{I_1}\!\!\cdot\cdot \cdot 
G_{I_{k-2}} {A}_{I_{k-1}} {W}_{I_{k}}\|
\leq 
\| g^{I_1}\|_{\infty}\!\!\cdot\cdot\cdot
\|g^{I_{k-2}}\|_{\infty}\,
 \|A_{I_{k-1}}{W}_{I_{k}}\|\leq  p_{r-s}^{\Cal R} \|{W}\|_{s}
+p_{r-s-1}^{\Cal R} \|{W}\|_{s+1}\tag 16.17
$$
By (4.9) applied to (16.10) in place of (4.3):
$$
|\dot{E}_I|\leq \big(1+\|\dot{g}\|_{\infty}+\|\pa \dot{p}\|_\infty/c_0\big) E_I
+2\sqrt{E_I} \|H_I\|\tag 16.18
$$
where $c_0$ is the constant in (1.6). By (16.14)-(16.17) we have 
$$
\|H_I\|\leq C
\sum_{s=0}^{r-1}
\big( \dot{m}_{r-s}^{\Cal R}\|\dot{W}\|_s + p_{r-s}^{\Cal R}  \|{W}\|_{s}\big)
+p_0^{\Cal R} \|{W}\|_{r}
+ \|{F}\|_{r} \tag 16.19
$$
and using (16.13) 
$$
\|H_I\|\leq K_1 
\sum_{s=0}^{r-1} \big(\dot{m}_{r-s}^{\Cal R}+p_{r-s}^{\Cal
R}\big)\big(C^{\Cal{R}}_s+E_{s}^{\Cal S}\big)
 +K_1 p_0^{\Cal R} \big(C^{\Cal{R}}_r+E_{r}^{\Cal S}\big)+ \|{F}\|_{r} \tag 16.20
$$
Summing (16.18) over all $I\in {\Cal S}$ with $|I|=r$ and using (16.20)
we get
$$\multline 
\Big| \frac{d E_r^{\Cal S}}{dt}\Big|\leq K_1 
\big(1+\|\dot{g}\|_{\infty}+\|\pa
\dot{p}\|_\infty/c_0+\sum_{S\in {\Cal S}}\|\pa Sp\|_\infty\big)
\big(C^{\Cal{R}}_r+E_{r}^{\Cal S}\big)
\\
+K_1 \sum_{s=0}^{r-1} \big(\dot{m}_{r-s}^{\Cal R}+p_{r-s}^{\Cal
R}\big)\big(C^{\Cal{R}}_s+E_{s}^{\Cal S}\big)+
\|{F}\|_{r} 
\endmultline\tag 16.21
$$
Furthermore, by Lemma {12.}1, (12.18) hold with ${\Cal U}$ replaced by 
${\Cal R}$ and ${\Cal T}$ replaced by ${\Cal S}$:
$$
\Big|\frac{d C^{\Cal{R}}_r}{dt}\Big|\leq K_1\dot{m}_0^{\Cal R}
 \big(C^{\Cal{R}}_r+ E^{\Cal S}_r\big) 
+K_1 \sum_{s=1}^{r-1} \dot{m}^{\Cal R}_{r-s}
\big( C^{\Cal{R}}_s+ E^{\Cal S}_s\big)+\|F\|_r\tag 16.22
$$
(16.21) together with (16.22) gives us a bound for 
$C^{\Cal{R}}_r+E^{\Cal{S}}_r$ in terms of $C^{\Cal{R}}_s+E^{\Cal{S}}_s$ for $s<r$:
$$\multline 
C^{\Cal{R}}_r(t)+E^{\Cal{S}}_r(t)
\leq K_1 e^{K_1\int_0^t n\, d\tau} \big(C^{\Cal{R}}_r(0)+E^{\Cal{S}}_r(0)\big)\\
+K_1 e^{K_1\int_0^t n\, d\tau} 
\int_0^t \Big(  \sum_{s=1}^{r-1} \big(\dot{m}_{r-s}^{\Cal R}+p_{r-s}^{\Cal
R}\big)\big(C^{\Cal{R}}_s+E_{s}^{\Cal S}\big)+
\|{F}\|_{r}   \Big) \, d\tau 
\endmultline \tag 16.23
$$
where $n=1+\|\dot{g}\|_{\infty}+\|\pa
\dot{p}\|_\infty/c_0+\sum_{S\in {\Cal S}}\|\pa Sp\|_\infty
+\|\omega\|_{\infty}$. 
Since we already have proven the bound for $E_0^{\Cal S}=E_0$ in section 4,
(16.23) inductively gives a bound for $C^{\Cal{R}}_r+E^{\Cal{S}}_r$. 
Hence by (16.13) we obtain: 
\proclaim{Lemma {16.}1}  Suppose that
$x, p\in C^{r+2}([0,T]\times\Omega)$, $p\,\big|_{\pa\Omega}\!\!=0$,
$\na_N p\,\big|_{\pa\Omega}\!\leq \!-c_0\!<\!0$ and $\div V\!=\!0$, where $V\!=\!D_t x$.
Let $W$ be the solution of (16.1) where $F$ is divergence free.
Then there is a constant $C$ depending only on the norm
of $(x,p)$, a lower bound for the constant $c_0$ and an upper bound for $T$, 
such that, for $0\!\leq\! t\!\leq\! T$, we have
$$
\| \dot{W}(t)\|_{r}+\|{W}(t)\|_{r}+\langle W(t)\rangle_{A,r}
\leq C \Big(\| \dot{W}(0)\|_{r}+\|{W}(0)\|_{r}+\langle W(0)\rangle_{A,r}
 +\int_0^t \!\!  \|F\|_{r}\, d\tau \Big)
\tag 16.24
$$
where
$$
\|W(t)\|_r=\sum_{|I|\leq r,\, I\in{\Cal R}} \|{\Cal L}_U^I W(t)\|_{L^2(\Omega)},
\qquad
\langle W(t)\rangle_{A,r}=\sum_{|I|\leq r,\, I\in {\Cal S}}
 \langle {\Cal L}_S^I W(t),A{\Cal L}_S^I W(t)\rangle^{1/2}\tag 16.25
$$
\endproclaim

Note that the $\|W(t)\|_r$ is equivalent to the usual time independent
Sobolev norm. Since there are compactly supported divergence free vector fields
$\langle W(t)\rangle_{A,r}$ is only a semi-norm on divergence free vector fields, see (3.10).
Furthermore, since $0\! <\! c_0\!\leq\!  -\na_N p\! \leq\!  C$ it follows from (3.11)
that $\langle W(t)\rangle_{A,r}$ is equivalent to a time
independent semi-norm given by (3.11) with $f$ the distance function $d(y)$,
see (6.2). 
Since we only apply tangential vector fields, it also follows from (3.11) 
that, up to lower order terms that can be bounded by $\|W(t)\|_r$, 
it is equivalent to that the normal 
component of the vector field $W_N=N_a W^a$ is in $H^r(\pa\Omega)$. 

\demo{Definition {16.}1} 
With notation as in (16.25) 
define $H^r(\Omega)$ to be the completion of $C^\infty(\Omega)$  in 
the norm $\|W(t)\|_r$
and define $N^r(\Omega)$ to be the completion of  the divergence free $C^\infty(\Omega)$ 
vector fields 
in the norm $\|W\|_{N^r}=\|W(t)\|_r+\langle W(t)\rangle_{A,r}$. 
\enddemo 

Since the projection onto divergence free vector fields is continuous in the $H^r$ norm 
it follows that $H^r$ is also the completion of the divergence free $C^\infty$ vector fields in 
the $H^r$ norm.

\proclaim{Theorem {16.}2}  Suppose that 
$x, p\in C^{r+2}([0,T]\times\Omega)$, $p\big|_{\pa\Omega}=0$,
$\na_N p\big|_{\pa\Omega}\leq -c_0<0$ and $\div D_t x=0$. 
Then if initial data and the inhomogeneous term in (2.29) are divergence
free and satisfy
$$
(W_0,W_1)\in  N^r(\Omega) \times  H^r(\Omega),\qquad
F\in L^1\big([0,T], H^r(\Omega)\big)
\tag 16.26
$$
the linearized equations (2.29) have a solution
$$
(W,\dot{W})\in C\big([0,T], N^r(\Omega)\times H^r(\Omega)\big)\tag 16.27
$$
\endproclaim 
\demo{Proof} The existence of a solution in (16.27) 
follows from section 15 if initial data and the inhomogeneous term are divergence free and  $C^\infty$ 
and the inhomogeneous term is supported in $t>0$. 
By approximating the initial data and the inhomogeneous term 
in (16.26) with $C^\infty$ divergence free vector fields 
and applying the estimate (16.24) to the differences we get a 
convergent sequence in (16.27) so the limit 
must also be in this space.
\qed\enddemo

\subheading{Acknowledgments} I would like to thank Salah Baouendi, Demetrios 
Christodoulou, David Ebin and Kate Okikiolu for helpful discussions.

\enddocument